\documentclass[twoside,reqno,A4]{amsart}%%%%%%%%%%%%%%%%%%%%%%
\usepackage{latexsym,rotate,eucal,cite,url}%%%%%%%%%%%%%%%%%%%
\usepackage{amsmath,amsthm,amssymb,amsxtra} %%%%{mtpro2} %%%%%

\newcommand{\quab}{\hspace*{2.0mm}}
\newcommand{\puab}{\hspace*{-1.5mm}}
\newcommand{\sst}{\scriptstyle}

\newcommand{\Lla}{\Longleftarrow}%

\newcommand{\pav}[1]{\lfloor{#1}\rfloor}

\newcommand{\pavv}[1]{\left\lfloor{#1}\right\rfloor}

\newcommand{\ang}[1]{\langle{#1}\rangle}

\newcommand{\mb}[1]{\mathbb{#1}}
\newcommand{\mc}[1]{\mathcal{#1}}

\newcommand{\hh}[2]{\left[#1\right]_{#2}}%%%{\left[#1;#2\right]_{#3}}
\newcommand{\hyp}[3]{\left[\puab\ba{#1}#2\\#3\ea\puab\right]}

\newcommand{\binm}{\binom}

\newcommand{\be}{\begin{equation}}
\newcommand{\ee}{\end{equation}}
\newcommand{\ba}{\begin{array}}
\newcommand{\ea}{\end{array}}
\newcommand{\bmn}{\begin{eqnarray}}
\newcommand{\emn}{\end{eqnarray}}
\newcommand{\bnm}{\begin{eqnarray*}}
\newcommand{\enm}{\end{eqnarray*}}
\newcommand{\bln}{\begin{subequations}}
\newcommand{\eln}{\end{subequations}}
\newcommand{\blm}{\bln\bmn}
\newcommand{\elm}{\emn\eln}

\newcommand{\pq}[1]{\begin{equation}#1\end{equation}}
\newcommand{\pp}[2]{\begin{aligned}[#1]#2 % [t]&[b] %%
            \end{aligned}}  %%%aligned in equation%%%%

\newcommand{\pnq}[1]{\begin{align*}#1%%%%%%%%%%%%%%%%%
            \end{align*}}    %%%align=auto-space%%%%%%
\newcommand{\pnp}[2]{\begin{alignat*}{#1}#2%%%%%%%%%%%
            \end{alignat*}}  %%%alignat=no-space%%%%%%
\newcommand{\centro}[1]
           {\begin{center}#1\end{center}}

\newcommand{\centerbox}[2]{\centro{\begin{tabular}{|c|}\hline
\parbox{#1}{\mbox{}\\[2.5mm]#2\mbox{}\\}\\\hline\end{tabular}}}

\newcommand{\alp}{\alpha}
\newcommand{\bet}{\beta}
\newcommand{\gam}{\gamma}

\newcommand{\lam}{\lambda}
\newcommand{\vph}{\varphi}

\newcommand{\Ome}{\Omega}

\newcommand{\Gam}{\Gamma}
\newtheorem{thm}{Theorem}%[section]
\newtheorem{lemm}[thm]{Lemma}

\newtheorem{exam}{Example}

\newtheorem{entry}{Entry}%%%%%%%%%%%%%%%%

\newcommand{\nota}[2]{\centerbox{#1}{\textbf{Achtung}\quad#2}}
	
\newcommand{\referxy}[4]{\bibitem{kn:#1}{#2,}~\emph{#3,}~{#4.}}	
\newcommand{\cito}[1]{\cite{kn:#1}}	
\newcommand{\citu}[2]{\cite[#2]{kn:#1}}

\begin{document}%{\fbox{\fns\today}\hfill} %%%%%% arXiv %%%%%%
%{\fns Lecce\ \copyright\ \textbf{University of Salento}}%%%%%
%{\fns Lecce\ \copyright\ \textbf{Working in Progress}}%%%%%%%
%%%%%%%%%%%%%%%%%%%%%%%%%%%%%%%%%%%%%%%%%%%%%%%%%%%%%%%%%%%%%%
\title{$\pi$-Formulae from Dual Series\\ of the Dougall Theorem}
\author{Wenchang Chu}
\thanks{Email address: chu.wenchang@unisalento.it}
\address{Department of Mathematics and Physics\newline
	University of Salento (P.~O.~Box~193) \newline
      	73100 Lecce, ~ Italy}
\email{chu.wenchang@unisalento.it}

%%%%%%%%%%%%%%%%%%%%%%%%%%%%%%%%%%%%%%%%%%%%%%%%%%%%%%%%%%%%%%
%%%%%%%%%%%%%%%%%%%%%%%%%%%%%%%%%%%%%%%%%%%%%%%%%%%%%%%%%%%%%%

\begin{abstract}
By means of the extended Gould-Hsu inverse series relations,
we find that the dual relation of Dougall's summation theorem
for the well--poised $_7F_6$-series can be utilized to construct
numerous interesting Ramanujan--like infinite series expressions 
for $\pi^{\pm1}$ and $\pi^{\pm2}$, including an elegant formula 
of $\pi^{-2}$ due to Guillera.
\end{abstract}

%%%%%%%%%%%%%%%%%%%%%%%%%%%%%%%%%%%%%%%%%%%%%%%%%%%%%%%%%%%%%%%%%%%%%%%
%%%%%%%%%%%%%%%%%%%%%%%%%%%%%%%%%%%%%%%%%%%%%%%%%%%%%%%%%%%%%%%%%%%%%%%
\maketitle\thispagestyle{empty}%%%%%%%%\setcounter{page}{0}%%%%%%%%%%%%
\markboth{Wenchang Chu}{Dual Series of the Dougall Theorem} %%%%%%%%
%%%%%%%%%%%%%%%%%%%%%%%%%%%%%%%%%%%%%%%%%%%%%%%%%%%%%%%%%%%%%%%%%%%%%%%
%%%%%%%%%%%%%%%%%%%%%%%%%%%%%%%%%%%%%%%%%%%%%%%%%%%%%%%%%%%%%%%%%%%%%%%

%%%%%%%%%%%%%%%%%%%%%%%%%%%%%%%%%%%%%%%%%%%%%%%%%%%%%%%%%%%%%%%%%%%%%%%
\vspace*{-5mm}\section{Introduction and Motivation}
%%%%%%%%%%%%%%%%%%%%%%%%%%%%%%%%%%%%%%%%%%%%%%%%%%%%%%%%%%%%%%%%%%%%%%%
In 1973, Gould and Hsu~\cito{hsu} discovered a useful pair
of inverse series relations, which can equivalently be
reproduced below. Let $\{a_i,b_i\}$ be any two complex
sequences such that the $\vph$-polynomials defined by
\pq{\vph(x;0)\equiv1\label{phi-pp}
\quad\text{and}\quad
\vph(x;n)=\prod_{k=0}^{n-1}(a_k+xb_k)
\quad\text{for}\quad n\in\mb{N}}
differ from zero for $x,\:n\in\mb{N}_0$.
Then there hold the inverse series relations
\blm
f(n)&=&\label{inv+gh}
\sum_{k=0}^{n}
(-1)^k\binm{n}{k}
\vph(k;n)\:g(k),\\
g(n)&=&\label{inv-gh}
\sum_{k=0}^{n}
(-1)^k\binm{n}{k}
\frac{a_k+kb_k}{\vph(n;k+1)}
\:f(k).
\elm
This inverse pair has wide applications to terminating
hypergeometric series identities
\cite{kn:chu93b,kn:chu93c,kn:chu94a,kn:chu94b,kn:chu94e,kn:gess82s}.
The duplicate form with applications can be found
in \cite{kn:chu00a,kn:chu02b,kn:chu08j}. There
exist also $q$-analogues due to Carlitz~\cito{carlitz}
which has applications to $q$-series identities
\cite{kn:chu94c,kn:chu94d,kn:chu95b,kn:chu06d,kn:gess83s,kn:gess86s}.

The Gould--Hsu inversions have the following extended form
(cf.~\cite{kn:bress,kn:chu93b,kn:chu94e}):
\blm
f(n)&=&\label{inv+w}
\sum_{k=0}^{n}(-1)^k\binm{n}{k}
\vph(\lam+k;n)\vph(-k;n)
\frac{\lam+2k}{(\lam+n)_{k+1}}g(k),\\
g(n)&=&\label{inv-w}
\sum_{k=0}^{n}
(-1)^k\binm{n}{k}
\frac{(a_k+\lam b_k+kb_k)(a_k-kb_k)}{\vph(\lam+n;k+1)\vph(-n;k+1)}(\lam+k)_n
\:f(k);
\elm
where the shifted factorials are defined by
\[(x)_0=1
\quad\text{and}\quad
(x)_n=\frac{\Gam(x+n)}{\Gam(x)}
=x(x+1)\cdots(x+n-1)
\quad\text{for}\quad
n\in\mb{N}.\]
%%%%%%%%%%%%%%%%%%%%%%%%%%%%%%%%%%%%%%%%%%%%%%%%%%%%%%%%%%%%%%%%%%%%%%%
There exist numerous hypergeometric series identities
(see~\citu{bry}{Chapter~8}%%%~\citu{table}{Chapter~7}
and \cite{kn:chu93b,kn:chu93c,kn:chu94a,kn:chu94b,kn:chu94e,
          kn:chu17k,kn:chu18g,kn:ges95,kn:gess82s} for example).
One of well--known summation theorems originally due to Dougall~\cito{dougall}
is about the terminating well--poised $_7F_6$-series. By examining its dual formulae
through (\ref{inv+w}--\ref{inv-w}), we find that their limiting relations result
unexpectedly in $\pi$-related infinite series expressions, including 
the following elegant formula discovered
%%%%%%%%%%%%%%%%%%  by Ramanujan~\cito{ramanujan} %%%%%%%%%%%%%%%%%%%%%
%\pnq{\frac{4}{\pi}
%=&\sum_{k=0}^{\infty}
%\hyp{ccc}
%{\frac12,\frac12,\frac12}{\rule[1mm]{0mm}{3mm}1,\:1,\:1}_k
%\frac{1+6k}{4^k},\\
%\frac{8}{\pi}
%=&\sum_{k=0}^{\infty}
%\hyp{ccc}
%{\frac12,\frac14,\frac34}
%{\rule[1mm]{0mm}{3mm}1,\:1,\:1}_k
%\frac{3+20k}{(-4)^k},\\
%\frac{16}{\pi}\!
%=&\sum_{k=0}^{\infty}
%\hyp{ccc}
%{\frac12,\frac12,\frac12}{\rule[1mm]{0mm}{3mm}1,\:1,\:1}_k
%\frac{5+42k}{64^k},}
by Guillera~\cite{kn:gj03em,kn:gj06rj,kn:gj08rj}
%%%%%%%%%%%%%%%%%%%%%%%%%%%%%%%%%%%%%%%%%%%%%%%%%%%%%%%%%%%%%%%%%%%%%%%
\[\frac{32}{\pi^2}
=\sum_{k=0}^{\infty}
\hyp{ccccccc}{\frac12,\:\frac12,\:\frac12,\:\frac14,\:\frac34}
{1,\:1,\:1,\:1,\:1\rule[2mm]{0mm}{2mm}}_k
\frac{3+34k+120k^2}{16^k}.\]%%%\xquad\fbox{\eqref{pi+32}}
%\[\frac{3 \pi ^2}{32}=1+\sum_{k=1}^{\infty}\hyp{c}
%{1,-\frac{1}{2},-\frac{1}{2},-\frac{1}{4},-\frac{3}{4}\\[-3mm]}
%{\frac{5}{2},~\frac{3}{4},~\frac{3}{4},~\frac{5}{4},~\frac{5}{4}}_k
%\frac{3+3k-22 k^2-40 k^3}{16^k}.\qquad\fbox{\eqref{pi-32}}\]
%%%%%%%%%%%%%%%%%%%%%%%%%%%%%%%%%%%%%%%%%%%%%%%%%%%%%%%%%%%%%%%%%%%%%%%
By means of the duplicate forms of (\ref{inv+w}--\ref{inv-w}), we shall work out,
in details, the dual formulae of Dougall's summation theorem in the next section.
Then applications will be presented in Section~3, where several $\pi$-related 
infinite series of Ramanujan--like~\cito{ramanujan} with the convergence rate 
``$\frac1{16}$" will be illustrated as examples.
%%%%%%%%%%%%%%%%%%%%%%%%%%%%%%%%%%%%%%%%%%%%%%%%%%%%%%%%%%%%%%%%%%%%%%%
%Further results will be illustrated by carrying out the same procedure to the
%triplicate, quadruplicate and quintuplicate forms of (\ref{inv+w}--\ref{inv-w}).

Recall that the $\Gam$-function~(see~\citu{rain}{\S8}
for example) is defined by the beta integral
\[\Gam(x)=\int_{0}^{\infty}u^{x-1}
e^{-u}\mathrm{d}u
\quad\text{for}\quad
\Re(x)>0,\]
which admits Euler's reflection property
\pq{\Gam(x)\Gam(1-x)=\label{reflex}
\frac{\pi}{\sin\pi x}
\quad\text{with}\quad \Gam(\tfrac12)=\sqrt{\pi}.}
The following asymptotic formula
\pq{\Gam(x+n)\approx n^{x}(n-1)!
\quad\text{as}\quad\label{limit}
n\to\infty,}
will be useful in evaluating limits of $\Gam$-function quotients.

For the sake of brevity, the product and quotient
of shifted factorials will respectively be abbreviated to
\pnq{[\alp,\bet,\cdots,\gam]_n\quab
=&(\alp)_n(\bet)_n\cdots(\gam)_n,\\
\hyp{cccc}{\alp,\bet,\cdots,\gam}
          {A,B,\cdots,C}_n
=&\frac{(\alp)_n(\bet)_n\cdots(\gam)_n}
        {(A)_n(B)_n\cdots(C)_n}.}
The similar notation will be employed for the $\Gam$-function quotient
\[\Gam\hyp{cccc}{\alp,\bet,\cdots,\gam}
          {A,B,\cdots,C}
=\frac{\Gam(\alp)\Gam(\bet)\cdots\Gam(\gam)}
      {\Gam(A)\Gam(B)\cdots\Gam(C)}.\]

\section{Main Theorems from Duplicate Inversions}
%%%%%%%%%%%%%%%%%%%%%%%%%%%%%%%%%%%%%%%%%%%%%%%%%%%%%%%%%%%%%%%%%%%%%%%
%%%%%%%%%%%%%%%%%%%%%%%%%%%%%%%%%%%%%%%%%%%%%%%%%%%%%%%%%%%%%%%%%%%%%%%

%%%\[\boxed{n=\pavv{\tfrac{n}2}+\pavv{\tfrac{1+n}2}}\]
The fundamental identity discovered by Dougall~\cito{dougall}
(see also Bailey~\citu{bailey}{\S4.3}) for very--well--poised
terminating $_7F_6$-series can be stated as
\pq{\label{doug}
\pp{c}{\Ome_n(a;b,c,d)
~:=~\hyp{c}{1+a,1+a-b-c,1+a-b-d,1+a-c-d}{1+a-b,1+a-c,1+a-d,1+a-b-c-d}_n&\\
=\sum_{k=0}^n\frac{a+2k}{a}\hyp{ccccc}
{a,~b,~c,~d,~e,\:-n}
{1,1+a-b,1+a-c,1+a-d,1+a-e,1+a+n}_k&,}}
where the series is 2-balanced because $1+2a+n=b+c+d+e$.

For all $n\in\mb{N}_0$, it is well known that
$n=\pavv{\tfrac{n}2}+\pavv{\tfrac{1+n}2}$, where $\pav{x}$
denotes the greatest integer not exceeding $x$.
%%%%%%%%%%%%%%%%%%%%%%%%%%%%%%%%%%%%%%%%%%%%%%%%%%%%%%%%%%%%%%%%%%%%%%%%%%%
%%%%%%%%%%%%%%%%%%%%%%%%%%%%%%%%%%%%%%%%%%%%%%%%%%%%%%%%%%%%%%%%%%%%%%%%%%%
Then it is not difficult to check that Dougall's formula \eqref{doug}
is equivalent to the following one
\pnq{\Ome_n\big(a;b+\pavv{\tfrac{n}2},c,d+\pavv{\tfrac{1+n}2}\big)
=&\hyp{c}{1+a-c-d,b+c-a}{1+a-d,b-a}_{\pavv{\frac{n}2}}\\
\times\hyp{c}{1+a,\quad b+d-a}{1+a-c,b+c+d-a}_n
&\hyp{c}{1+a-b-c,c+d-a}{1+a-b,d-a}_{\pavv{\frac{1+n}2}}}
with its parameters subject to $\boxed{1+2a=b+c+d+e}$.
Reformulate the above equality as a binomial sum
\pnq{\sum_{k=0}^n&(-1)^k\binm{n}{k}
\hh{b+k,b-a-k}{\pavv{\frac{n}2}}
\hh{d+k,d-a-k}{\pavv{\frac{1+n}2}}\\
&~\times\frac{a+2k}{(a+n)_{k+1}}
\hyp{c}{a,~b,~c,~d,~1+2a-b-c-d}{1+a-b,1+a-c,1+a-d,b+c+d-a}_k\\
=~&\hyp{c}{b,1+a-c-d,b+c-a}{1+a-d}_{\pavv{\frac{n}2}}
~\!\hyp{c}{d,1+a-b-c,c+d-a}{1+a-b}_{\pavv{\frac{1+n}2}}\\
&~\times\hyp{c}{a,~b+d-a}{1+a-c,b+c+d-a}_n.}

This equality matches exactly to \eqref{inv+w} under the assignments
$\lam\to a$ and
\[\vph(x;n)=(b-a+x)_{\pavv{\frac{n}2}}(d-a+x)_{\pavv{\frac{1+n}2}}\]
as well as
\pnq{
f(n)=&\hyp{c}{1+a-c-d,b,b+c-a}{1+a-d}_{\pavv{\frac{n}2}}
\hyp{c}{a,~b+d-a}{1+a-c,b+c+d-a}_n\\
\times&\hyp{c}{1+a-b-c,d,c+d-a}{1+a-b}_{\pavv{\frac{1+n}2}},\\
g(k)=&\hyp{c}{a,~b,~c,~d,~1+2a-b-c-d}{1+a-b,1+a-c,1+a-d,b+c+d-a}_k.}
The dual relation corresponding to \eqref{inv-w} can explicitly
be stated, according to the parity of $k$ and
$(a)_k(a+k)_n=(a)_n(a+n)_k$, as
\pnq{&\hyp{c}
{b,~c,~d,~1+2a-b-c-d}
{1+a-b,1+a-c,1+a-d,b+c+d-a}_n\\
&=\sum_{k\ge0}\binm{n}{2k}
\frac{(d+3k)(d-a-k)(a+n)_{2k}}
{\hh{b+n,b-a-n}{k}\hh{d+n,d-a-n}{k+1}}\\
&\times\hyp{c}{1+a-c-d,b,b+c-a}{1+a-d}_{k}
\hyp{c}{1+a-b-c,d,c+d-a}{1+a-b}_{k}\\
&\times\hyp{c}{b+d-a}{1+a-c,b+c+d-a}_{2k}\\
&-\sum_{k\ge0}\binm{n}{2k+1}
\frac{(b+3k+1)(b-a-k-1)(a+n)_{2k+1}}
{\hh{b+n,b-a-n}{k+1}\hh{d+n,d-a-n}{k+1}}\\
&\times\hyp{c}{1+a-c-d,b,b+c-a}{1+a-d}_{k}
\hyp{c}{1+a-b-c,d,c+d-a}{1+a-b}_{k+1}\\
&\times\hyp{c}{b+d-a}{1+a-c,b+c+d-a}_{2k+1}.}

Now multiplying by $``n^2"$ across this binomial relation
and then letting $n\to\infty$, we may evaluate the limits
of the left member by \eqref{limit} and of the corresponding 
right member through the Weierstrass $M$-test on uniformly 
convergent series (cf. Stromberg~\citu{karl}{\S3.106}).
After some routine simplification, the resulting limiting
relation can be expressed explicitly in the following lemma.
%%%%%%%%%%%%%%%%%%%%%%%%%%%%%%%%%%%%%%%%%%%%%%%%%%%%%%%%%%%%%%%%%%%%%%%%%%%
\begin{lemm}[Infinite series identity]\label{2=1+1}
\pnq{\Gam&\hyp{c}
{1+a-b,1+a-c,1+a-d,b+c+d-a}
{b,~c,~d,~1+2a-b-c-d}\\
&=\sum_{k\ge0}\frac{(d+3k)(a-d)}{(2k)!}
\hyp{c}{b+d-a}{1+a-c,b+c+d-a}_{2k}\\
&\times\hyp{c}{1+a-c-d,b,b+c-a}{a-d}_{k}
\hyp{c}{1+a-b-c,d,c+d-a}{1+a-b}_{k}\\
&+\sum_{k\ge0}\frac{(b+3k+1)(a-b)}{(2k+1)!}
\hyp{c}{b+d-a}{1+a-c,b+c+d-a}_{2k+1}\\
&\times\hyp{c}{1+a-c-d,b,b+c-a}{1+a-d}_{k}
\hyp{c}{1+a-b-c,d,c+d-a}{a-b}_{k+1}.}
\end{lemm}

According to this lemma, we are going to show two main theorems
that will be utilized, in the next section, to deduce infinite
series expressions for $\pi^{\pm1}$ and $\pi^{\pm2}$.

For the equality in Lemma~\ref{2=1+1}, multiplying both
sides by $(1+a-c)(b+c+d-a)$ and then unifying the two
sums, we derive the following infinite series identity.
%%%%%%%%%%%%%%%%%%%%%%%%%%%%%%%%%%%%%%%%%%%%%%%%%%%%%%%%%%%%%%%%%%%%%%%%%%%
\begin{thm}[Infinite series identity]\label{2=1+1a}
\pnq{\Gam&
\hyp{c}{1+a-b,2+a-c,1+a-d,1-a+b+c+d}
{b,\quad c,\quad d,\quad1+2a-b-c-d}\\
=&\sum_{k=0}^{\infty}\mc{P}(k)
\frac{[b,d,1+a-b-c,1+a-c-d,b+c-a,c+d-a]_k(b+d-a)_{2k}}
{(2k+1)!~[1+a-b,1+a-d]_k[2+a-c,1-a+b+c+d]_{2k}},}
where $\mc{P}(k)$ is the polynomial given by
\pnq{\mc{P}(k)
&=(1+a-b-c+k)(d+k)(c+d-a+k)(b+d-a+2k)(1+b+3k)\\
&+(1+2k)(a-d+k)(1+a-c+2k)(b+c+d-a+2k)(d+3k).}
\end{thm}

Alternatively, by shifting backward $k\to k-1$ for the second
sum and then unifying it to the first one, we get analogously,
from Lemma~\ref{2=1+1} another infinite series identity.
%%%%%%%%%%%%%%%%%%%%%%%%%%%%%%%%%%%%%%%%%%%%%%%%%%%%%%%%%%%%%%%%%%%%%%%%%%%
\begin{thm}[Infinite series identity]\label{2=1+1b}
\pnq{\Gam&\label{2=1+1}
\hyp{c}{1+a-b,1+a-c,1+a-d,b+c+d-a}
{b,\quad c,\quad d,\quad1+2a-b-c-d}\\
=&\sum_{k=0}^{\infty}\mc{Q}(k)
\frac{[b,d,1+a-b-c,1+a-c-d,b+c-a,c+d-a]_k(b+d-a)_{2k}}
{(2k)!~[1+a-b,1+a-d]_k[1+a-c,b+c+d-a]_{2k}},}
where $\mc{Q}(k)$ is the rational function given by
\[\mc{Q}(k)=(a-d+k)(d+3k)
\bigg\{1+\tfrac
{(2k)(a-b+k)(a-c+2k)(b+c+d-a-1+2k)(b-2+3k)}
{\sst(a-c-d+k)(b-1+k)(b+c-a-1+k)(b+d-a-1+2k)(d+3k)}\bigg\}.\]
\end{thm}

%\nota{111mm}{When running with \emph{Mathematica},
%the initial condition $\mc{Q}(0)=d(a-d)$ should be imposed,
%because the fraction inside the braces may differ from zero
%for specific values of $\{a,b,c,d\}$.}

%%%%%%%%%%%%%%%%%%%%%%%%%%%%%%%%%%%%%%%%%%%%%%%%%%%%%%%%%%%%%%%%%%%%%%%%%%%
%%%%%%%%%%%%%%%%%%%%%%%%%%%%%%%%%%%%%%%%%%%%%%%%%%%%%%%%%%%%%%%%%%%%%%%%%%%
%%%%%%%%%%%%%%%%%%%%%%%%%%%%%%%%%%%%%%%%%%%%%%%%%%%%%%%%%%%%%%%%%%%%%%%%%%%
%%%%%%%%%%%%%%%%%%%%%%%%%%%%%%%%%%%%%%%%%%%%%%%%%%%%%%%%%%%%%%%%%%%%%%%%%%%
%%%%% \input{check.tex} %%%%%%%%%%%%%%%%%%%%%%%%%%%%%%%%%%%%%%%%%%%%%%%%%%%
%%%%%%%%%%%%%%%%%%%%%%%%%%%%%%%%%%%%%%%%%%%%%%%%%%%%%%%%%%%%%%%%%%%%%%%%%%%
%%%%%%%%%%%%%%%%%%%%%%%%%%%%%%%%%%%%%%%%%%%%%%%%%%%%%%%%%%%%%%%%%%%%%%%%%%%
%%%%%%%%%%%%%%%%%%%%%%%%%%%%%%%%%%%%%%%%%%%%%%%%%%%%%%%%%%%%%%%%%%%%%%%%%%%
%%%%%%%%%%%%%%%%%%%%%%%%%%%%%%%%%%%%%%%%%%%%%%%%%%%%%%%%%%%%%%%%%%%%%%%%%%%
\section{Infinite Series for $\pi^{\pm1}$ and $\pi^{\pm2}$}
%%%%%%%%%%%%%%%%%%%%%%%%%%%%%%%%%%%%%%%%%%%%%%%%%%%%%%%%%%%%%%%%%%%%%%%%%%%
%%%%%%%%%%%%%%%%%%%%%%%%%%%%%%%%%%%%%%%%%%%%%%%%%%%%%%%%%%%%%%%%%%%%%%%%%%%
%%%%%%%%%%%%%%%%%%%%%%%%%%%%%%%%%%%%%%%%%%%%%%%%%%%%%%%%%%%%%%%%%%%%%%%%%%%
By applying Theorems~\ref{2=1+1a} and~\ref{2=1+1b}, we can derive numerous
infinite series identities. They are recorded below in seven classes whose 
weight polynomial degrees are not greater than $3$. 
For all the examples, the parameter settings $\boxed{a,b,c,d}$ and eventual 
references are highlight in their headers. In order to ensure the accuracy,
all the summation formulae in this section are verified experimentally 
by appropriately devised \emph{Mathematica} commands.

\subsection{Series for $\pi^{-2}$}
%%%%%%%%%%%%%%%%%%%%%%%%%%%%%%%%%%%%%%%%%%%%%%%%%%%%%%%%%%%%%%%%%%%%%%%%%%%
%%%%%%%%%%%%%%%%%%%%%%%%%%%%%%%%%%%%%%%%%%%%%%%%%%%%%%%%%%%%%%%%%%%%%%%%%%%
%%%%%%%%%%%%%%%%%%%%%%%%%%%%%%%%%%%%%%%%%%%%%%%%%%%%%%%%%%%%%%%%%%%%%%%%%%%

\begin{exam}[Guillera~\cite{kn:gj03em,kn:gj06rj,kn:gj08rj}:
\fbox{$\frac12,\frac12,\frac12,\frac12$} in Theorem~\ref{2=1+1a}]\label{pi+32}
\[ \frac{32}{\pi ^2}
=\sum_{k=0}^{\infty}
\hyp{c}
{ \frac{1}{2},\frac{1}{2},\frac{1}{2},\frac{1}{4},\frac{3}{4}}
{\\[-3mm] 1,1,1,1,1}_k
\frac{120 k^2+34 k+3}{16^k}.\]
\end{exam}

%%%%%%%%%%%%%%%%%%%%%%%%%%%%%%%%%%%%%%%%%%%%%%%%%%%%%%%%%%%%%%%%%%%%%%%%%%%
\begin{exam}[Chu and Zhang~\cito{chu14mc}:
\fbox{$\frac32,\frac12,\frac32,\frac12$} in Theorem~\ref{2=1+1a}]
\[ \frac{128}{\pi ^2}
=\sum_{k=0}^{\infty}
\hyp{c}
{ \frac{1}{2},\frac{1}{2},\frac{1}{2},\frac{1}{4},-\frac{1}{4}}
{\\[-3mm] 1,1,1,~2,~2}_k
\frac{120 k^2+118 k+13}{16^k}.\]
\end{exam}

%%%%%%%%%%%%%%%%%%%%%%%%%%%%%%%%%%%%%%%%%%%%%%%%%%%%%%%%%%%%%%%%%%%%%%%%%%%
\begin{exam}[\fbox{$\frac32,\frac12,\frac12,\frac32$} in Theorem~\ref{2=1+1a}]
\[ \frac{256}{3 \pi ^2}
=\sum_{k=0}^{\infty}
\hyp{c}
{\frac{1}{2},-\frac{1}{2},\frac{3}{2},\frac{1}{4},\frac{3}{4}}
{\\[-3mm] 1,1,1,~2,~2}_k
\frac{80 k^3+148 k^2+80 k+9}{16^k}.\]
\end{exam}

%%%%%%%%%%%%%%%%%%%%%%%%%%%%%%%%%%%%%%%%%%%%%%%%%%%%%%%%%%%%%%%%%%%%%%%%%%%
\begin{exam}[\fbox{$\frac32,\frac32,\frac12,\frac32$} in Theorem~\ref{2=1+1a}]
\[ \frac{512}{\pi ^2}
=\sum_{k=0}^{\infty}
\hyp{c}
{ \frac{1}{2},\frac{1}{2},\frac{3}{2},\frac{3}{4},\frac{5}{4}}
{\\[-3mm] 1,1,1,2,2}_k
\frac{240 k^3+532 k^2+336 k+45}{16^k}.\]
\end{exam}

%%%%%%%%%%%%%%%%%%%%%%%%%%%%%%%%%%%%%%%%%%%%%%%%%%%%%%%%%%%%%%%%%%%%%%%%%%%
\begin{exam}[\fbox{$\frac12,\frac12,\frac12,-\frac12$} in Theorem~\ref{2=1+1a}]
\[ \frac{32}{\pi ^2}
=\sum_{k=0}^{\infty}
\hyp{c}
{\frac{3}{2},-\frac{1}{2},-\frac{1}{2},\frac{1}{4},-\frac{1}{4},\frac{7}{6}}
{\\[-3mm] 1,~ 1,~ 1,~1,~ 2,~\frac{1}{6}}_k
\frac{3-10 k-40 k^2}{16^k}.\]
\end{exam}

%%%%%%%%%%%%%%%%%%%%%%%%%%%%%%%%%%%%%%%%%%%%%%%%%%%%%%%%%%%%%%%%%%%%%%%%%%%
\begin{exam}[\fbox{$\frac32,\frac12,\frac32,-\frac12$} in Theorem~\ref{2=1+1a}]
\[ \frac{256}{3 \pi ^2}
=\sum_{k=0}^{\infty}
\hyp{c}
{\frac{3}{2},-\frac{1}{2},-\frac{1}{2},-\frac{1}{4},-\frac{3}{4},\frac{7}{6}}
{\\[-3mm] 1,~ 1,~ 1,~ 2,~3,~\frac{1}{6}}_k
\frac{9-38 k-40 k^2}{16^k}.\]
\end{exam}

%%%%%%%%%%%%%%%%%%%%%%%%%%%%%%%%%%%%%%%%%%%%%%%%%%%%%%%%%%%%%%%%%%%%%%%%%%%
\begin{exam}[\fbox{$\frac32,\frac32,\frac12,\frac32$} in Theorem~\ref{2=1+1b}]
\[\frac{8}{\pi ^2}
=\sum_{k=1}^{\infty}
\hyp{c}
{\frac{3}{2},-\frac{1}{2},-\frac{1}{2},\frac{1}{4},\frac{3}{4},\frac{7}{6}}
{\\[-3mm]1,~ 1,~ 1,~ 1,~ 1,~ \frac{1}{6}}_k
\frac{k(3-18 k+40 k^2)}{16^k}.\]
\end{exam}

%%%%%%%%%%%%%%%%%%%%%%%%%%%%%%%%%%%%%%%%%%%%%%%%%%%%%%%%%%%%%%%%%%%%%%%%%%%
\begin{exam}[\fbox{$\frac32,\frac12,\frac32,\frac12$} in Theorem~\ref{2=1+1b}]
\[\frac{24}{\pi ^2}
=\sum_{k=0}^{\infty}
\hyp{c}
{\frac{1}{2},-\frac{1}{2},-\frac{1}{2},-\frac{1}{4},-\frac{3}{4},\frac{5}{4},\frac{5}{6}}
{\\[-3mm] 1,~ 1,~ 1,~ 1,~ 2,\frac{1}{4},-\frac{1}{6}}_k
\frac{3+8 k+20 k^2}{16^k}.\]
\end{exam}

%%%%%%%%%%%%%%%%%%%%%%%%%%%%%%%%%%%%%%%%%%%%%%%%%%%%%%%%%%%%%%%%%%%%%%%%%%%
\begin{exam}[\fbox{$\frac32,\frac32,-\frac12,\frac12$} in Theorem~\ref{2=1+1a}]
\[ \frac{256}{9 \pi ^2}
=\sum_{k=0}^{\infty}
\hyp{c}
{\frac{3}{2},\frac{5}{2},-\frac{1}{2},-\frac{3}{2},\frac{1}{4},\frac{3}{4}}
{\\[-3mm] 1,~ 1,~ 1,~ 2,~2,~\frac{1}{2}}_k
\frac{5+12 k-68 k^2-80 k^3}{16^k}.\]
\end{exam}

\subsection{Series for $\pi^2$} \ %% \fbox{bisection series?}
%%%%%%%%%%%%%%%%%%%%%%%%%%%%%%%%%%%%%%%%%%%%%%%%%%%%%%%%%%%%%%%%%%%%%%%%%%%
%%%%%%%%%%%%%%%%%%%%%%%%%%%%%%%%%%%%%%%%%%%%%%%%%%%%%%%%%%%%%%%%%%%%%%%%%%%
%%%%%%%%%%%%%%%%%%%%%%%%%%%%%%%%%%%%%%%%%%%%%%%%%%%%%%%%%%%%%%%%%%%%%%%%%%%

\begin{exam}[Chu and Zhang~\cito{chu14mc}:
\fbox{$\frac32,1,1,1$} in Theorem~\ref{2=1+1a}]
\[ \frac{9 \pi ^2}{8}
=\sum_{k=0}^{\infty}
\hyp{c}
{ 1,\frac{1}{2},\frac{1}{2},\frac{1}{4},\frac{3}{4}}
{\\[-3mm] \frac{3}{2},\frac{5}{4},\frac{5}{4},\frac{7}{4},\frac{7}{4}}_k
\frac{11+64 k+111 k^2+60 k^3}{16^k}.\]
\end{exam}

%%%%%%%%%%%%%%%%%%%%%%%%%%%%%%%%%%%%%%%%%%%%%%%%%%%%%%%%%%%%%%%%%%%%%%%%%%%
\begin{exam}[\fbox{$\frac52,2,1,2$} in Theorem~\ref{2=1+1a}]
\[ \frac{225 \pi ^2}{32}
=\sum_{k=0}^{\infty}
\hyp{c}
{ 2,2,\frac{1}{2},\frac{1}{2},\frac{3}{4},\frac{5}{4}}
{\\[-3mm] 1,\frac{3}{2},\frac{7}{4},\frac{7}{4},\frac{9}{4},\frac{9}{4}}_k
\frac{68+206 k+197 k^2+60 k^3 }{16^k}.\]
\end{exam}

%%%%%%%%%%%%%%%%%%%%%%%%%%%%%%%%%%%%%%%%%%%%%%%%%%%%%%%%%%%%%%%%%%%%%%%%%%%
\begin{exam}[\fbox{$\frac52,1,2,2$} in Theorem~\ref{2=1+1a}]
\[ \frac{135 \pi ^2}{64}
=\sum_{k=0}^{\infty}
\hyp{c}
{ 2,\frac{1}{2},-\frac{1}{2},\frac{5}{3},\frac{1}{4},\frac{3}{4}}
{\\[-3mm] \frac{5}{2},\frac{2}{3},~\frac{5}{4},\frac{7}{4},\frac{7}{4},\frac{9}{4}}_k
\frac{21+93 k+110 k^2+40 k^3 }{16^k}.\]
\end{exam}

%%%%%%%%%%%%%%%%%%%%%%%%%%%%%%%%%%%%%%%%%%%%%%%%%%%%%%%%%%%%%%%%%%%%%%%%%%%
\begin{exam}[\fbox{$\frac52,1,2,1$} in Theorem~\ref{2=1+1b}]\label{pi-32}
\[\frac{3 \pi ^2}{32}
=1+\sum_{k=1}^{\infty}
\hyp{c}
{1,-\frac{1}{2},-\frac{1}{2},-\frac{1}{4},-\frac{3}{4}}
{\\[-3mm]\frac{5}{2},~\frac{3}{4},~\frac{3}{4},~\frac{5}{4},~\frac{5}{4}}_k
\frac{3+3k-22 k^2-40 k^3}{16^k}.\]
\end{exam}

%%%%%%%%%%%%%%%%%%%%%%%%%%%%%%%%%%%%%%%%%%%%%%%%%%%%%%%%%%%%%%%%%%%%%%%%%%%
\begin{exam}[\fbox{$\frac72,1,2,1$} in Theorem~\ref{2=1+1b}]
\[\frac{15 \pi ^2}{256}
=\frac13+\sum_{k=1}^{\infty}
\hyp{c}
{1,-\frac{1}{2},-\frac{3}{2},-\frac{3}{4},-\frac{5}{4}}
{\\[-3mm] \frac{7}{2},~ \frac{1}{4},~ \frac{3}{4},~ \frac{5}{4},~ \frac{7}{4}}_k
\frac{1-3 k+2 k^2+8 k^3}{16^k}.\]
\end{exam}

%%%%%%%%%%%%%%%%%%%%%%%%%%%%%%%%%%%%%%%%%%%%%%%%%%%%%%%%%%%%%%%%%%%%%%%%%%%
\begin{exam}[\fbox{$\frac72,2,2,2$} in Theorem~\ref{2=1+1b}]
\[\frac{27 \pi ^2}{128}
=\sum_{k=0}^{\infty}
\hyp{c}
{2,-\frac{1}{2},-\frac{1}{2},\frac{4}{3},\frac{1}{4},-\frac{1}{4}}
{\\[-3mm] \frac{5}{2},~ \frac{1}{3},~ \frac{5}{4},~ \frac{5}{4},~ \frac{7}{4},~ \frac{7}{4}}_k
\frac{2-21 k-66 k^2-40 k^3}{16^k}.\]
\end{exam}

%%%%%%%%%%%%%%%%%%%%%%%%%%%%%%%%%%%%%%%%%%%%%%%%%%%%%%%%%%%%%%%%%%%%%%%%%%%
\begin{exam}[\fbox{$\frac72,1,2,2$} in Theorem~\ref{2=1+1b}]
\[\frac{405 \pi ^2}{256}
=18+\sum_{k=1}^{\infty}
\hyp{c}
{2,-\frac{1}{2},\frac{3}{2},-\frac{3}{2},-\frac{1}{4},-\frac{3}{4}}
{\\[-3mm] \frac{1}{2},~ \frac{7}{2},~ \frac{3}{4},~ \frac{5}{4},~ \frac{5}{4},~ \frac{7}{4}}_k
\frac{48-59 k-194 k^2-120 k^3}{16^k}.\]
\end{exam}

\subsection{Series for $\pi^2/\Gam^3$}
%%%%%%%%%%%%%%%%%%%%%%%%%%%%%%%%%%%%%%%%%%%%%%%%%%%%%%%%%%%%%%%%%%%%%%%%%%%
%%%%%%%%%%%%%%%%%%%%%%%%%%%%%%%%%%%%%%%%%%%%%%%%%%%%%%%%%%%%%%%%%%%%%%%%%%%
%%%%%%%%%%%%%%%%%%%%%%%%%%%%%%%%%%%%%%%%%%%%%%%%%%%%%%%%%%%%%%%%%%%%%%%%%%%

\begin{exam}[\fbox{$\frac12,\frac13,\frac13,-\frac23$} in Theorem~\ref{2=1+1a}]
\[ \frac{98 \pi ^2}{3 \Gamma (\frac{2}{3})^3}
=\sum_{k=0}^{\infty}
\hyp{c}
{\frac{1}{3},-\frac{2}{3},\frac{5}{6},-\frac{5}{6},
\frac{11}{6},\frac{10}{9},\frac{1}{12},-\frac{5}{12}}
{\\[-3mm] 1,~ \frac{3}{2},~\frac{1}{4},~\frac{3}{4},~
\frac{13}{6},~\frac{1}{9},~\frac{13}{12},~\frac{19}{12}}_k
\frac{118+45 k-1098 k^2-1080 k^3 }{16^k}.\]
\end{exam}

%%%%%%%%%%%%%%%%%%%%%%%%%%%%%%%%%%%%%%%%%%%%%%%%%%%%%%%%%%%%%%%%%%%%%%%%%%%
\begin{exam}[\fbox{$\frac32,\frac13,\frac13,\frac43$} in Theorem~\ref{2=1+1a}]
\[ \frac{637 \pi ^2}{16 \Gamma (\frac{2}{3})^3}
=\sum_{k=0}^{\infty}
\hyp{c}
{ \frac{1}{3},\frac{4}{3},\frac{5}{6},-\frac{5}{6},
\frac{11}{6},\frac{13}{9},\frac{1}{12},\frac{7}{12}}
{\\[-3mm] 1,\frac{3}{2},\frac{3}{4},~\frac{5}{4},
\frac{13}{6},\frac{4}{9},\frac{19}{12},\frac{25}{12}}_k
\frac{1080 k^3+2286 k^2+1395 k+161}{16^k}.\]
\end{exam}

%%%%%%%%%%%%%%%%%%%%%%%%%%%%%%%%%%%%%%%%%%%%%%%%%%%%%%%%%%%%%%%%%%%%%%%%%%%
\begin{exam}[\fbox{$\frac12,\frac23,-\frac13,-\frac13$} in Theorem~\ref{2=1+1a}]
\[\frac{275 \pi ^2}{\Gamma (\frac{1}{3})^3}
=\sum_{k=0}^{\infty}
\hyp{c}
{\frac{2}{3},-\frac{1}{3},\frac{7}{6},-\frac{7}{6},
\frac{13}{6},\frac{11}{9},-\frac{1}{12},\frac{5}{12}}
{\\[-3mm] 1,~\frac{3}{2},~\frac{1}{4},~\frac{3}{4},~
\frac{11}{6},~\frac{2}{9},~\frac{17}{12},~\frac{23}{12}}_k
\frac{125-351 k-1602 k^2-1080 k^3}{16^k}.\]
\end{exam}

%%%%%%%%%%%%%%%%%%%%%%%%%%%%%%%%%%%%%%%%%%%%%%%%%%%%%%%%%%%%%%%%%%%%%%%%%%%
\begin{exam}[\fbox{$\frac32,\frac23,\frac53,-\frac13$} in Theorem~\ref{2=1+1a}]
\[ \frac{825 \pi ^2}{8 \Gamma (\frac{1}{3})^3}
=\sum_{k=0}^{\infty}
\hyp{c}
{ \frac{2}{3},-\frac{1}{3},\frac{1}{6},-\frac{1}{6},
\frac{7}{6},\frac{11}{9},-\frac{1}{12},-\frac{7}{12}}
{\\[-3mm] 1,~ \frac{3}{2},~\frac{3}{4},~\frac{5}{4},~
\frac{17}{6},~\frac{2}{9},~\frac{11}{12},~\frac{17}{12}}_k
\frac{53-315 k-1278 k^2-1080 k^3}{16^k}.\]
\end{exam}

%%%%%%%%%%%%%%%%%%%%%%%%%%%%%%%%%%%%%%%%%%%%%%%%%%%%%%%%%%%%%%%%%%%%%%%%%%%
\begin{exam}[\fbox{$\frac32,\frac23,-\frac13,\frac53$} in Theorem~\ref{2=1+1a}]
\[ \frac{2805 \pi ^2}{4 \Gamma (\frac{1}{3})^3}
=\sum_{k=0}^{\infty}
\hyp{c}
{ \frac{2}{3},\frac{5}{3},\frac{7}{6},-\frac{7}{6},
\frac{13}{6},\frac{14}{9},\frac{5}{12},\frac{11}{12}}
{\\[-3mm] 1,\frac{3}{2},\frac{3}{4},~\frac{5}{4},
\frac{11}{6},\frac{5}{9},\frac{23}{12},\frac{29}{12}}_k
\frac{1080 k^3+2790 k^2+2151 k+478}{16^k}.\]
\end{exam}

%%%%%%%%%%%%%%%%%%%%%%%%%%%%%%%%%%%%%%%%%%%%%%%%%%%%%%%%%%%%%%%%%%%%%%%%%%%
\begin{exam}[\fbox{$-\frac12,-\frac23,-\frac23,-\frac23$} in Theorem~\ref{2=1+1b}]
\[\frac{3872 \pi ^2}{243 \Gamma (\frac{2}{3})^3}
=\sum_{k=0}^{\infty}
\hyp{c}
{ -\frac{2}{3},-\frac{5}{3},\frac{5}{6},\frac{11}{6},
-\frac{11}{6},\frac{4}{9},-\frac{5}{12},-\frac{11}{12}}
{\\[-3mm] 1,~\frac{1}{2},-\frac{1}{4},-\frac{3}{4},~
\frac{7}{6},-\frac{5}{9},~\frac{7}{12},~\frac{13}{12} }_k
\frac{1080 k^3-954 k^2-585 k+242}{16^k}.\]
\end{exam}

%%%%%%%%%%%%%%%%%%%%%%%%%%%%%%%%%%%%%%%%%%%%%%%%%%%%%%%%%%%%%%%%%%%%%%%%%%%
\begin{exam}[\fbox{$\frac12,-\frac23,\frac43,-\frac23$} in Theorem~\ref{2=1+1b}]
\[ \frac{2380 \pi ^2}{27 \Gamma (\frac{2}{3})^3}
=\sum_{k=0}^{\infty}
\hyp{c}
{-\frac{2}{3},-\frac{5}{3},-\frac{1}{6},-\frac{5}{6},
\frac{5}{6},\frac{4}{9},-\frac{11}{12},-\frac{17}{12}}
{\\[-3mm] 1,~ \frac{1}{2},~\frac{1}{4},-\frac{1}{4},~
\frac{13}{6},-\frac{5}{9},~\frac{1}{12},~\frac{7}{12}}_k
\frac{1080 k^3-1278 k^2+99 k+170}{16^k}.\]
\end{exam}

%%%%%%%%%%%%%%%%%%%%%%%%%%%%%%%%%%%%%%%%%%%%%%%%%%%%%%%%%%%%%%%%%%%%%%%%%%%
\begin{exam}[\fbox{$\frac12,\frac43,-\frac23,-\frac23$} in Theorem~\ref{2=1+1b}]
\[ \frac{770 \pi ^2}{27 \Gamma (\frac{2}{3})^3}
=\sum_{k=0}^{\infty}
\hyp{c}
{\frac{1}{3},-\frac{2}{3},\frac{5}{6},\frac{11}{6},
 -\frac{11}{6},\frac{7}{9},-\frac{5}{12},\frac{1}{12}}
{\\[-3mm] 1,~ \frac{1}{2},~\frac{1}{4},-\frac{1}{4},
\frac{7}{6},-\frac{2}{9},~\frac{13}{12},~\frac{19}{12}}_k
\frac{55+441 k-234 k^2-1080 k^3}{16^k}.\]
\end{exam}

%%%%%%%%%%%%%%%%%%%%%%%%%%%%%%%%%%%%%%%%%%%%%%%%%%%%%%%%%%%%%%%%%%%%%%%%%%%
\begin{exam}[\fbox{$\frac32,\frac43,-\frac23,\frac43$} in Theorem~\ref{2=1+1b}]
\[ \frac{-1001 \pi ^2}{18 \Gamma (\frac{2}{3})^3}
=\sum_{k=0}^{\infty}
\hyp{c}
{ \frac{1}{3},\frac{4}{3},\frac{5}{6},\frac{11}{6},
-\frac{11}{6},\frac{10}{9},\frac{1}{12},\frac{7}{12}}
{\\[-3mm] 1,\frac{1}{2},~\frac{1}{4},~\frac{3}{4},~
\frac{7}{6},~\frac{1}{9},~\frac{19}{12},\frac{25}{12} }_k
\frac{1080 k^3+1422 k^2+351 k+44}{16^k}.\]
\end{exam}

%%%%%%%%%%%%%%%%%%%%%%%%%%%%%%%%%%%%%%%%%%%%%%%%%%%%%%%%%%%%%%%%%%%%%%%%%%%
\begin{exam}[\fbox{$\frac32,\frac43,\frac43,-\frac23$} in Theorem~\ref{2=1+1b}]
\[ \frac{385 \pi ^2}{36 \Gamma (\frac{2}{3})^3}
=\sum_{k=0}^{\infty}
\hyp{c}
{\frac{1}{3},-\frac{2}{3},-\frac{1}{6},-\frac{5}{6},
\frac{5}{6},\frac{7}{9},-\frac{5}{12},-\frac{11}{12} }
{\\[-3mm] 1,~ \frac{1}{2},~\frac{1}{4},~\frac{3}{4},~\frac{13}{6},
-\frac{2}{9},~\frac{7}{12},~\frac{13}{12}}_k
\frac{55+189 k+90 k^2-1080 k^3}{16^k}.\]
\end{exam}

%%%%%%%%%%%%%%%%%%%%%%%%%%%%%%%%%%%%%%%%%%%%%%%%%%%%%%%%%%%%%%%%%%%%%%%%%%%
\begin{exam}[\fbox{$\frac12,-\frac13,\frac23,-\frac13$} in Theorem~\ref{2=1+1b}]
\[ \frac{910\pi ^{2}}{9\Gamma(\frac{1}{3})^3}
=\sum_{k=0}^{\infty}
\hyp{c}
{ -\frac{1}{3},-\frac{4}{3},\frac{1}{6},\frac{7}{6},
-\frac{7}{6},\frac{5}{9},-\frac{7}{12},-\frac{13}{12}}
{\\[-3mm] 1,~ \frac{1}{2},~\frac{1}{4},-\frac{1}{4},~
\frac{11}{6},-\frac{4}{9},~\frac{5}{12},~\frac{11}{12} }_k
\frac{1080 k^3-774 k^2-225 k+91}{16^k}.\]
\end{exam}

%%%%%%%%%%%%%%%%%%%%%%%%%%%%%%%%%%%%%%%%%%%%%%%%%%%%%%%%%%%%%%%%%%%%%%%%%%%
\begin{exam}[\fbox{$\frac32,\frac23,\frac23,\frac23$} in Theorem~\ref{2=1+1b}]
\[ \frac{1225 \pi ^2}{6 \Gamma (\frac{1}{3})^3}
=\sum_{k=0}^{\infty}
\hyp{c}
{ \frac{2}{3},-\frac{1}{3},\frac{1}{6},\frac{7}{6},
-\frac{7}{6},\frac{8}{9},-\frac{1}{12},-\frac{7}{12}}
{\\[-3mm] 1,\frac{1}{2},~\frac{1}{4},~\frac{3}{4},~
\frac{11}{6},-\frac{1}{9},~\frac{11}{12},~\frac{17}{12}}_k
\frac{98+153 k-414 k^2-1080 k^3}{16^k}.\]
\end{exam}

\subsection{Series for $\Gam^3/\pi^2$}
%%%%%%%%%%%%%%%%%%%%%%%%%%%%%%%%%%%%%%%%%%%%%%%%%%%%%%%%%%%%%%%%%%%%%%%%%%%
%%%%%%%%%%%%%%%%%%%%%%%%%%%%%%%%%%%%%%%%%%%%%%%%%%%%%%%%%%%%%%%%%%%%%%%%%%%
%%%%%%%%%%%%%%%%%%%%%%%%%%%%%%%%%%%%%%%%%%%%%%%%%%%%%%%%%%%%%%%%%%%%%%%%%%%

\begin{exam}[\fbox{$-\frac12,-\frac56,\frac16,\frac16$} in Theorem~\ref{2=1+1a}]
\[\frac{180 \Gamma (\frac{2}{3})^3}{\pi ^2}
=\sum_{k=0}^{\infty}
\hyp{c}
{ \frac{1}{6},-\frac{1}{6},\frac{5}{6},-\frac{5}{6},
\frac{7}{6},-\frac{1}{12},\frac{5}{12},\frac{19}{18}}
{\\[-3mm] 1,~ 1,~ \frac{1}{2},~\frac{3}{2},~\frac{1}{3},~
\frac{2}{3},~\frac{4}{3},~\frac{1}{18}}_k
\frac{35+228 k-540 k^2-2160 k^3}{16^k}.\]
\end{exam}

%%%%%%%%%%%%%%%%%%%%%%%%%%%%%%%%%%%%%%%%%%%%%%%%%%%%%%%%%%%%%%%%%%%%%%%%%%%
\begin{exam}[\fbox{$-\frac12,\frac16,\frac16,\frac76$} in Theorem~\ref{2=1+1a}]
\[ \frac{8748 \Gamma (\frac{2}{3})^3}{7 \pi ^2}
=\sum_{k=0}^{\infty}
\hyp{c}
{\frac{1}{6},\frac{5}{6},-\frac{5}{6},\frac{7}{6},
\frac{11}{6},\frac{11}{12},\frac{17}{12},\frac{25}{18}}
{\\[-3mm] 1,2,\frac{3}{2},~\frac{3}{2},~\frac{1}{3},~
\frac{2}{3},-\frac{2}{3},~ \frac{7}{18}}_k
\frac{593+1344 k-1404 k^2-2160 k^3}{16^k}.\]
\end{exam}

%%%%%%%%%%%%%%%%%%%%%%%%%%%%%%%%%%%%%%%%%%%%%%%%%%%%%%%%%%%%%%%%%%%%%%%%%%%
\begin{exam}[\fbox{$\frac12,\frac16,-\frac56,\frac76$} in Theorem~\ref{2=1+1a}]
\[ \frac{960 \Gamma (\frac{2}{3})^3}{7 \pi ^2}
=\sum_{k=0}^{\infty}
\hyp{c}
{ \frac{1}{6},-\frac{1}{6},\frac{7}{6},-\frac{7}{6},
\frac{13}{6},\frac{5}{12},\frac{11}{12},\frac{25}{18}}
{\\[-3mm] 1,~ 1,~ \frac{1}{2},~\frac{3}{2},~\frac{1}{3},~
\frac{4}{3},~\frac{5}{3},~\frac{7}{18}}_k
\frac{65+372 k-756 k^2-2160 k^3}{16^k}.\]
\end{exam}

%%%%%%%%%%%%%%%%%%%%%%%%%%%%%%%%%%%%%%%%%%%%%%%%%%%%%%%%%%%%%%%%%%%%%%%%%%%
\begin{exam}[\fbox{$\frac12,\frac16,\frac76,-\frac56$} in Theorem~\ref{2=1+1a}]
\[ \frac{960 \Gamma (\frac{2}{3})^3}{\pi ^2}
=\sum_{k=0}^{\infty}
\hyp{c}
{ \frac{1}{6},-\frac{1}{6},\frac{5}{6},-\frac{5}{6},
 \frac{7}{6},-\frac{1}{12},-\frac{7}{12},\frac{19}{18}}
{\\[-3mm] 1,~ 1,~ \frac{1}{2},~\frac{3}{2},~\frac{2}{3},~
\frac{4}{3},~\frac{7}{3},~\frac{1}{18} }_k
\frac{245-204 k-2052 k^2-2160 k^3}{16^k}.\]
\end{exam}

%%%%%%%%%%%%%%%%%%%%%%%%%%%%%%%%%%%%%%%%%%%%%%%%%%%%%%%%%%%%%%%%%%%%%%%%%%%
\begin{exam}[\fbox{$\frac12,\frac16,\frac76,\frac76$} in Theorem~\ref{2=1+1a}]
\[ \frac{7776 \Gamma (\frac{2}{3})^3}{7 \pi ^2}
=\sum_{k=0}^{\infty}
\hyp{c}
{\frac{1}{6},\frac{5}{6},\frac{5}{6},-\frac{5}{6},\frac{7}{6},
   \frac{11}{6},\frac{5}{12},\frac{11}{12},\frac{25}{18}}
{\\[-3mm] 1,2,~\frac{3}{2},\frac{3}{2},~\frac{1}{3},
\frac{2}{3},\frac{4}{3},-\frac{1}{6},~\frac{7}{18}}_k
\frac{360 k^2+546 k+191}{16^k}.\]
\end{exam}

%%%%%%%%%%%%%%%%%%%%%%%%%%%%%%%%%%%%%%%%%%%%%%%%%%%%%%%%%%%%%%%%%%%%%%%%%%%
\begin{exam}[\fbox{$\frac32,\frac16,\frac16,\frac76$} in Theorem~\ref{2=1+1a}]
\[ \frac{1024 \Gamma (\frac{2}{3})^3}{21 \pi ^2}
=\sum_{k=0}^{\infty}
\hyp{c}
{\frac{1}{6},-\frac{1}{6},\frac{7}{6},-\frac{7}{6},
\frac{13}{6},-\frac{1}{12},\frac{5}{12},\frac{25}{18}}
{\\[-3mm] 1,~ 1,~ \frac{1}{2},~\frac{3}{2},~\frac{4}{3},~
\frac{5}{3},~\frac{7}{3},~\frac{7}{18}}_k
\frac{2160 k^3+2268 k^2+60 k+13}{16^k}.\]
\end{exam}

%%%%%%%%%%%%%%%%%%%%%%%%%%%%%%%%%%%%%%%%%%%%%%%%%%%%%%%%%%%%%%%%%%%%%%%%%%%
\begin{exam}[\fbox{$-\frac12,-\frac16,\frac56,\frac56$} in Theorem~\ref{2=1+1a}]
\[ \frac{2916 \Gamma (\frac{1}{3})^3}{5 \pi ^2}
=\sum_{k=0}^{\infty}
\hyp{c}
{-\frac{1}{6},\frac{5}{6},\frac{7}{6},-\frac{7}{6},
\frac{13}{6},\frac{7}{12},\frac{13}{12},\frac{23}{18}}
{\\[-3mm] 1,2,~\frac{3}{2},\frac{3}{2},-\frac{1}{3},
\frac{1}{3},\frac{2}{3},~\frac{5}{18}}_k
\frac{697+1056 k-1836 k^2-2160 k^3}{16^k}.\]
\end{exam}

%%%%%%%%%%%%%%%%%%%%%%%%%%%%%%%%%%%%%%%%%%%%%%%%%%%%%%%%%%%%%%%%%%%%%%%%%%%
\begin{exam}[\fbox{$\frac12,-\frac16,\frac{11}6,\frac56$} in Theorem~\ref{2=1+1a}]
\[ \frac{2592 \Gamma (\frac{1}{3})^3}{25 \pi ^2}
=\sum_{k=0}^{\infty}
\hyp{c}
{-\frac{1}{6},\frac{5}{6},\frac{7}{6},-\frac{7}{6},
\frac{13}{6},\frac{1}{12},\frac{7}{12},\frac{23}{18}}
{\\[-3mm] 1,~2,~\frac{3}{2},~\frac{3}{2},~\frac{1}{3},
\frac{2}{3},\frac{5}{3},~\frac{5}{18}}_k
\frac{223-888 k-3348 k^2-2160 k^3}{16^k}.\]
\end{exam}

%%%%%%%%%%%%%%%%%%%%%%%%%%%%%%%%%%%%%%%%%%%%%%%%%%%%%%%%%%%%%%%%%%%%%%%%%%%
\begin{exam}[\fbox{$\frac12,\frac56,-\frac16,-\frac16$} in Theorem~\ref{2=1+1a}]
\[ \frac{32 \Gamma (\frac{1}{3})^3}{5 \pi ^2}
=\sum_{k=0}^{\infty}
\hyp{c}
{\frac{1}{6},-\frac{1}{6},\frac{5}{6},-\frac{5}{6},
\frac{11}{6},\frac{1}{12},\frac{7}{12},\frac{23}{18}}
{\\[-3mm] 1,~ 1,~ \frac{1}{2},~\frac{3}{2},~\frac{2}{3},~
\frac{4}{3},~\frac{5}{3},~\frac{5}{18}}_k
\frac{2160 k^3+1188 k^2-84 k+11}{16^k}.\]
\end{exam}

%%%%%%%%%%%%%%%%%%%%%%%%%%%%%%%%%%%%%%%%%%%%%%%%%%%%%%%%%%%%%%%%%%%%%%%%%%%
\begin{exam}[\fbox{$\frac12,\frac56,-\frac16,\frac{11}6$} in Theorem~\ref{2=1+1a}]
\[ \frac{2592 \Gamma (\frac{1}{3})^3}{55 \pi ^2}
=\sum_{k=0}^{\infty}
\hyp{c}
{ \frac{1}{6},-\frac{1}{6},\frac{5}{6},\frac{7}{6},
\frac{11}{6},\frac{13}{12},\frac{19}{12},\frac{29}{18}}
{\\[-3mm] 1,2,~\frac{3}{2},\frac{3}{2},~\frac{2}{3},
\frac{4}{3},-\frac{1}{3},~\frac{11}{18}}_k
\frac{151+264 k-2052 k^2-2160 k^3}{16^k}.\]
\end{exam}

%%%%%%%%%%%%%%%%%%%%%%%%%%%%%%%%%%%%%%%%%%%%%%%%%%%%%%%%%%%%%%%%%%%%%%%%%%%
\begin{exam}[\fbox{$\frac32,-\frac16,\frac56,\frac56$} in Theorem~\ref{2=1+1a}]
\[\frac{256 \Gamma (\frac{1}{3})^3}{9 \pi ^2}
=\sum_{k=0}^{\infty}
\hyp{c}
{\frac{1}{6},-\frac{1}{6},\frac{5}{6},-\frac{5}{6},
\frac{11}{6},\frac{1}{12},-\frac{5}{12},\frac{23}{18}}
{\\[-3mm] 1,~ 1,~ \frac{1}{2},~\frac{3}{2},~\frac{4}{3},~
\frac{5}{3},~\frac{8}{3},~\frac{5}{18}}_k
\frac{55-348 k-2700 k^2-2160 k^3}{16^k}.\]
\end{exam}

%%%%%%%%%%%%%%%%%%%%%%%%%%%%%%%%%%%%%%%%%%%%%%%%%%%%%%%%%%%%%%%%%%%%%%%%%%%
\begin{exam}[\fbox{$\frac32,\frac56,\frac56,\frac{11}6$} in Theorem~\ref{2=1+1a}]
\[ \frac{6912 \Gamma (\frac{1}{3})^3}{55 \pi ^2}
=\sum_{k=0}^{\infty}
\hyp{c}
{\frac{1}{6},-\frac{1}{6},\frac{5}{6},\frac{7}{6},
\frac{11}{6},\frac{7}{12},\frac{13}{12},\frac{29}{18}}
{\\[-3mm] 1,2,~\frac{3}{2},\frac{3}{2},~\frac{2}{3},
\frac{4}{3},\frac{5}{3},~\frac{11}{18}}_k
\frac{2160 k^3+3564 k^2+1680 k+251}{16^k}.\]
\end{exam}

%%%%%%%%%%%%%%%%%%%%%%%%%%%%%%%%%%%%%%%%%%%%%%%%%%%%%%%%%%%%%%%%%%%%%%%%%%%
\begin{exam}[\fbox{$-\frac12,\frac16,\frac76,\frac16$} in Theorem~\ref{2=1+1b}]
\[ \frac{2673 \Gamma (\frac{2}{3})^3}{16 \pi ^2}
=\sum_{k=0}^{\infty}
\hyp{c}
{\frac{1}{6},\frac{5}{6},-\frac{5}{6},\frac{11}{6},
-\frac{11}{6},-\frac{1}{12},\frac{5}{12},\frac{13}{18}}
{\\[-3mm] 1,~ 1,~ \frac{1}{2},~\frac{3}{2},~\frac{1}{3},
-\frac{1}{3},-\frac{2}{3},-\frac{5}{18}}_k
\frac{11+1380 k+1188 k^2-2160 k^3}{16^k}.\]
\end{exam}

%%%%%%%%%%%%%%%%%%%%%%%%%%%%%%%%%%%%%%%%%%%%%%%%%%%%%%%%%%%%%%%%%%%%%%%%%%%
\begin{exam}[\fbox{$-\frac12,\frac76,-\frac56,\frac76$} in Theorem~\ref{2=1+1b}]
\[ \frac{13365 \Gamma (\frac{2}{3})^3}{16 \pi ^2}
=\sum_{k=0}^{\infty}
\hyp{c}
{\frac{1}{6},-\frac{1}{6},\frac{5}{6},-\frac{5}{6},
\frac{7}{6},\frac{11}{12},\frac{17}{12},\frac{19}{18}}
{\\[-3mm] 1,1,~ \frac{1}{2},\frac{3}{2},~  \frac{2}{3},
-\frac{2}{3},-\frac{5}{3},\frac{1}{18}}_k
\frac{2160 k^3-2484 k^2-1092 k+385}{16^k}.\]
\end{exam}

%%%%%%%%%%%%%%%%%%%%%%%%%%%%%%%%%%%%%%%%%%%%%%%%%%%%%%%%%%%%%%%%%%%%%%%%%%%
\begin{exam}[\fbox{$\frac12,\frac76,\frac16,\frac76$} in Theorem~\ref{2=1+1b}]
\[ \frac{675 \Gamma (\frac{2}{3})^3}{2 \pi ^2}
=\sum_{k=0}^{\infty}
\hyp{c}
{\frac{1}{6},-\frac{1}{6},\frac{5}{6},-\frac{5}{6},
\frac{7}{6},\frac{5}{12},\frac{11}{12},\frac{19}{18}}
{\\[-3mm] 1,1,~ \frac{1}{2},\frac{3}{2},~ \frac{1}{3},
\frac{2}{3},-\frac{2}{3},~ \frac{1}{18}}_k
\frac{2160 k^3-972 k^2-660 k+175}{16^k}.\]
\end{exam}

%%%%%%%%%%%%%%%%%%%%%%%%%%%%%%%%%%%%%%%%%%%%%%%%%%%%%%%%%%%%%%%%%%%%%%%%%%%
\begin{exam}[\fbox{$\frac32,\frac76,\frac76,\frac76$} in Theorem~\ref{2=1+1b}]
\[\frac{180 \Gamma (\frac{2}{3})^3}{\pi ^2}
=\sum_{k=0}^{\infty}
\hyp{c}
{\frac{1}{6},-\frac{1}{6},\frac{5}{6},-\frac{5}{6},
\frac{7}{6},-\frac{1}{12},\frac{5}{12},\frac{19}{18} }
{\\[-3mm] 1,~ 1,~ \frac{1}{2},~\frac{3}{2},~\frac{1}{3},~
\frac{2}{3},~ \frac{4}{3},~\frac{1}{18}}_k
\frac{35+228 k-540 k^2-2160 k^3}{16^k}.\]
\end{exam}

%%%%%%%%%%%%%%%%%%%%%%%%%%%%%%%%%%%%%%%%%%%%%%%%%%%%%%%%%%%%%%%%%%%%%%%%%%%
\begin{exam}[\fbox{$-\frac12,-\frac16,\frac{11}6,-\frac16$} in Theorem~\ref{2=1+1b}]
\[ \frac{1053 \Gamma (\frac{1}{3})^3}{32 \pi ^2}
=\sum_{k=0}^{\infty}
\hyp{c}
{-\frac{1}{6},\frac{7}{6},-\frac{7}{6},\frac{13}{6},
-\frac{13}{6},-\frac{5}{12},\frac{1}{12},\frac{11}{18}}
{\\[-3mm] 1,~ 1,~ \frac{1}{2},~\frac{3}{2},-\frac{1}{3},
-\frac{2}{3},~ \frac{2}{3},-\frac{7}{18}}_k
\frac{2160 k^3-756 k^2-1668 k+65}{16^k}.\]
\end{exam}

%%%%%%%%%%%%%%%%%%%%%%%%%%%%%%%%%%%%%%%%%%%%%%%%%%%%%%%%%%%%%%%%%%%%%%%%%%%
\begin{exam}[\fbox{$-\frac12,\frac56,-\frac16,\frac56$} in Theorem~\ref{2=1+1b}]
\[ \frac{3969 \Gamma (\frac{1}{3})^3}{32 \pi ^2}
=\sum_{k=0}^{\infty}
\hyp{c}
{\frac{1}{6},-\frac{1}{6},\frac{5}{6},\frac{7}{6},
 -\frac{7}{6},\frac{7}{12},\frac{13}{12},\frac{17}{18}}
{\\[-3mm] 1,1,\frac{1}{2},\frac{3}{2},~\frac{1}{3},
-\frac{1}{3},-\frac{4}{3},-\frac{1}{18}}_k
\frac{2160 k^3-2052 k^2-1092 k+245}{16^k}.\]
\end{exam}

%%%%%%%%%%%%%%%%%%%%%%%%%%%%%%%%%%%%%%%%%%%%%%%%%%%%%%%%%%%%%%%%%%%%%%%%%%%
\begin{exam}[\fbox{$\frac12,\frac56,\frac56,\frac56$} in Theorem~\ref{2=1+1b}]
\[ \frac{63 \Gamma (\frac{1}{3})^3}{10 \pi ^2}
=\sum_{k=0}^{\infty}
\hyp{c}
{ \frac{1}{6},-\frac{1}{6},\frac{5}{6},\frac{7}{6},
-\frac{7}{6},\frac{1}{12},\frac{7}{12},\frac{17}{18}}
{\\[-3mm] 1,~ 1,~ \frac{1}{2},\frac{3}{2},-\frac{1}{3},
\frac{1}{3},\frac{2}{3},-\frac{1}{18}}_k
\frac{432 k^3-108 k^2-132 k+7}{16^k}.\]
\end{exam}

%%%%%%%%%%%%%%%%%%%%%%%%%%%%%%%%%%%%%%%%%%%%%%%%%%%%%%%%%%%%%%%%%%%%%%%%%%%
\begin{exam}[\fbox{$\frac32,\frac{11}6,-\frac16,\frac{11}6$} in Theorem~\ref{2=1+1b}]
\[ \frac{84 \Gamma (\frac{1}{3})^3}{5 \pi ^2}
=\sum_{k=0}^{\infty}
\hyp{c}
{\frac{1}{6},-\frac{1}{6},\frac{5}{6},-\frac{5}{6},
\frac{11}{6},\frac{7}{12},\frac{13}{12},\frac{23}{18}}
{\\[-3mm] 1,1,~ \frac{1}{2},~\frac{3}{2},-\frac{1}{3},~
\frac{2}{3},~\frac{4}{3},~\frac{5}{18}}_k
\frac{2160 k^3-324 k^2-516 k+77}{16^k}.\]
\end{exam}

%%%%%%%%%%%%%%%%%%%%%%%%%%%%%%%%%%%%%%%%%%%%%%%%%%%%%%%%%%%%%%%%%%%%%%%%%%%
\begin{exam}[\fbox{$\frac32,\frac{11}6,\frac{11}6,-\frac16$} in Theorem~\ref{2=1+1b}]
\[\frac{84 \Gamma (\frac{1}{3})^3}{\pi ^2}
=\sum_{k=0}^{\infty}
\hyp{c}
{\frac{1}{6},-\frac{1}{6},\frac{5}{6},\frac{7}{6},-\frac{7}{6},\frac{1}{12},-\frac{5}{12},\frac{17}{18}}
{\\[-3mm] 1,~ 1,~ \frac{1}{2},~\frac{3}{2},~\frac{1}{3},~ \frac{2}{3},~\frac{5}{3},-\frac{1}{18}}_k
\frac{175+228 k-972 k^2-2160 k^3}{16^k}.\]
\end{exam}

\subsection{Series for $\pi^{-1}$}
%%%%%%%%%%%%%%%%%%%%%%%%%%%%%%%%%%%%%%%%%%%%%%%%%%%%%%%%%%%%%%%%%%%%%%%%%%%
%%%%%%%%%%%%%%%%%%%%%%%%%%%%%%%%%%%%%%%%%%%%%%%%%%%%%%%%%%%%%%%%%%%%%%%%%%%
%%%%%%%%%%%%%%%%%%%%%%%%%%%%%%%%%%%%%%%%%%%%%%%%%%%%%%%%%%%%%%%%%%%%%%%%%%%

\begin{exam}[Chu and Zhang~\cito{chu14mc}:
\fbox{$\frac12,\frac12,\frac12,\frac13$} in Theorem~\ref{2=1+1a}]
\[ \frac{15 \sqrt{3}}{\pi }
=\sum_{k=0}^{\infty}
\hyp{c}
{ \frac{1}{2},~\frac{1}{3},\frac{1}{3},\frac{2}{3},\frac{2}{3}}
{\\[-3mm] 1,1,1,\frac{11}{12},\frac{17}{12}}_k
\frac{135 k^2+75 k+8}{16^k}.\]
\end{exam}

%%%%%%%%%%%%%%%%%%%%%%%%%%%%%%%%%%%%%%%%%%%%%%%%%%%%%%%%%%%%%%%%%%%%%%%%%%%
\begin{exam}[\fbox{$\frac12,\frac12,\frac12,\frac23$} in Theorem~\ref{2=1+1a}]
\[ \frac{21 \sqrt{3}}{\pi }
=\sum_{k=0}^{\infty}
\hyp{c}
{ \frac{1}{2},~\frac{1}{3},\frac{1}{3},\frac{2}{3},\frac{2}{3}}
{\\[-3mm] 1,1,1,\frac{13}{12},\frac{19}{12}}_k
\frac{810 k^3+684 k^2+141 k+10 }{16^k}.\]
\end{exam}

%%%%%%%%%%%%%%%%%%%%%%%%%%%%%%%%%%%%%%%%%%%%%%%%%%%%%%%%%%%%%%%%%%%%%%%%%%%
\begin{exam}[\fbox{$\frac12,\frac12,\frac12,\frac14$} in Theorem~\ref{2=1+1a}]
\[ \frac{48}{\pi }
=\sum_{k=0}^{\infty}
\hyp{c}
{ \frac{1}{2},\frac{1}{4},\frac{3}{4},\frac{1}{8},\frac{5}{8}}
{\\[-3mm] 1,1,1,\frac{7}{8},\frac{11}{8}}_k
\frac{480 k^2+212 k+15}{16^k}.\]
\end{exam}

%%%%%%%%%%%%%%%%%%%%%%%%%%%%%%%%%%%%%%%%%%%%%%%%%%%%%%%%%%%%%%%%%%%%%%%%%%%
\begin{exam}[\fbox{$\frac12,\frac12,\frac12,\frac34$} in Theorem~\ref{2=1+1a}]
\[ \frac{80}{3 \pi }
=\sum_{k=0}^{\infty}
\hyp{c}
{ \frac{1}{2},\frac{1}{4},\frac{3}{4},\frac{3}{8},\frac{7}{8}}
{\\[-3mm] 1,1,1,\frac{9}{8},\frac{13}{8}}_k
\frac{640 k^3+560 k^2+112 k+7}{16^k}.\]
\end{exam}

%%%%%%%%%%%%%%%%%%%%%%%%%%%%%%%%%%%%%%%%%%%%%%%%%%%%%%%%%%%%%%%%%%%%%%%%%%%
\begin{exam}[\fbox{$\frac12,\frac12,\frac16,\frac12$} in Theorem~\ref{2=1+1a}]
\[ \frac{256}{3 \pi \sqrt{3}}
=\sum_{k=0}^{\infty}
\hyp{c}
{ \frac{1}{2},\frac{1}{4},\frac{3}{4},\frac{1}{6},\frac{5}{6}}
{\\[-3mm] 1,1,1,\frac{4}{3},\frac{5}{3} }_k
\frac{720 k^3+804 k^2+236 k+15}{16^k}.\]
\end{exam}

%%%%%%%%%%%%%%%%%%%%%%%%%%%%%%%%%%%%%%%%%%%%%%%%%%%%%%%%%%%%%%%%%%%%%%%%%%%
\begin{exam}[\fbox{$\frac12,\frac12,\frac12,\frac16$} in Theorem~\ref{2=1+1a}]
\[ \frac{192}{\pi\sqrt{3}}
=\sum_{k=0}^{\infty}
\hyp{c}
{ \frac{1}{2},\frac{1}{6},\frac{1}{6},\frac{1}{12},\frac{7}{12}}
{\\[-3mm] 1,1,1,~ \frac{4}{3},\frac{4}{3}}_k
\frac{6480 k^3+4284 k^2+840 k+35}{16^k}.\]
\end{exam}

%%%%%%%%%%%%%%%%%%%%%%%%%%%%%%%%%%%%%%%%%%%%%%%%%%%%%%%%%%%%%%%%%%%%%%%%%%%
\begin{exam}[Chu and Zhang~\cito{chu14mc}:
\fbox{$\frac12,\frac12,\frac12,\frac56$} in Theorem~\ref{2=1+1a}]
\[ \frac{384}{\pi\sqrt{3} }
=\sum_{k=0}^{\infty}
\hyp{c}
{ \frac{1}{2},\frac{5}{6},\frac{5}{6},\frac{5}{12},\frac{11}{12}}
{\\[-3mm] 1,1,1,~ \frac{2}{3},\frac{5}{3}}_k
\frac{1080 k^2+798 k+55}{16^k}.\]
\end{exam}

%%%%%%%%%%%%%%%%%%%%%%%%%%%%%%%%%%%%%%%%%%%%%%%%%%%%%%%%%%%%%%%%%%%%%%%%%%%
\begin{exam}[\fbox{$\frac12,\frac12,\frac12,\frac15$} in Theorem~\ref{2=1+1a}]
\[ \frac{105 \sqrt{5-2\sqrt{5}}}{\pi}
=\sum_{k=0}^{\infty}
\hyp{c}
{ \frac{1}{2},\frac{1}{5},\frac{1}{5},\frac{3}{5},\frac{4}{5},\frac{1}{10}}
{\\[-3mm] 1,1,1,\frac{13}{10},\frac{17}{20},\frac{27}{20}}_k
\frac{3750 k^3+2525 k^2+505 k+24}{16^k}.\]
\end{exam}

%%%%%%%%%%%%%%%%%%%%%%%%%%%%%%%%%%%%%%%%%%%%%%%%%%%%%%%%%%%%%%%%%%%%%%%%%%%
\begin{exam}[\fbox{$\frac12,\frac12,\frac12,\frac25$} in Theorem~\ref{2=1+1a}]
\[ \frac{45 \sqrt{5+2\sqrt{5}}}{\pi}
=\sum_{k=0}^{\infty}
\hyp{c}
{ \frac{1}{2},\frac{1}{5},\frac{2}{5},\frac{2}{5},\frac{3}{5},~\frac{7}{10}}
{\\[-3mm] 1,1,1,\frac{11}{10},\frac{19}{20},\frac{29}{20}}_k
\frac{3750 k^3+2800 k^2+595 k+42}{16^k}.\]
\end{exam}

%%%%%%%%%%%%%%%%%%%%%%%%%%%%%%%%%%%%%%%%%%%%%%%%%%%%%%%%%%%%%%%%%%%%%%%%%%%
\begin{exam}[\fbox{$\frac12,\frac12,\frac12,\frac35$} in Theorem~\ref{2=1+1a}]
\[\frac{55 \sqrt{5+2\sqrt{5}}}{3\pi}
=\sum_{k=0}^{\infty}
\hyp{c}
{ \frac{1}{2},\frac{2}{5},\frac{3}{5},\frac{3}{5},\frac{4}{5},~\frac{3}{10}}
{\\[-3mm] 1,1,1,\frac{9}{10},\frac{21}{20},\frac{31}{20}}_k
\frac{1250 k^3+1025 k^2+215 k+16}{16^k}.\]
\end{exam}

%%%%%%%%%%%%%%%%%%%%%%%%%%%%%%%%%%%%%%%%%%%%%%%%%%%%%%%%%%%%%%%%%%%%%%%%%%%
\begin{exam}[\fbox{$\frac12,\frac12,\frac12,\frac45$} in Theorem~\ref{2=1+1a}]
\[ \frac{195 \sqrt{5-2\sqrt{5}}}{\pi}
=\sum_{k=0}^{\infty}
\hyp{c}
{ \frac{1}{2},\frac{1}{5},\frac{2}{5},\frac{4}{5},\frac{4}{5},~\frac{9}{10}}
{\\[-3mm] 1,1,1,\frac{7}{10},\frac{23}{20},\frac{33}{20}}_k
\frac{3750 k^3+3350 k^2+655 k+36}{16^k}.\]
\end{exam}

%%%%%%%%%%%%%%%%%%%%%%%%%%%%%%%%%%%%%%%%%%%%%%%%%%%%%%%%%%%%%%%%%%%%%%%%%%%
\begin{exam}[\fbox{$\frac12,\frac12,\frac12,\frac18$} in Theorem~\ref{2=1+1a}]
\[ \frac{480}{\pi (\sqrt2+1)}
=\sum_{k=0}^{\infty}
\hyp{c}
{ \frac{1}{2},\frac{1}{8},\frac{1}{8},\frac{7}{8},\frac{1}{16},\frac{9}{16}}
{\\[-3mm] 1,1,1,\frac{11}{8},\frac{13}{16},\frac{21}{16}}_k
\frac{15360 k^3+9920 k^2+1888 k+63}{16^k}.\]
\end{exam}

%%%%%%%%%%%%%%%%%%%%%%%%%%%%%%%%%%%%%%%%%%%%%%%%%%%%%%%%%%%%%%%%%%%%%%%%%%%
\begin{exam}[\fbox{$\frac12,\frac12,\frac12,\frac38$} in Theorem~\ref{2=1+1a}]
\[ \frac{224}{3 \pi(\sqrt2-1) }
=\sum_{k=0}^{\infty}
\hyp{c}
{ \frac{1}{2},\frac{3}{8},\frac{3}{8},\frac{5}{8},\frac{3}{16},\frac{11}{16}}
{\\[-3mm] 1,1,1,\frac{9}{8},\frac{15}{16},\frac{23}{16}}_k
\frac{5120 k^3+3776 k^2+800 k+55}{16^k}.\]
\end{exam}

%%%%%%%%%%%%%%%%%%%%%%%%%%%%%%%%%%%%%%%%%%%%%%%%%%%%%%%%%%%%%%%%%%%%%%%%%%%
\begin{exam}[\fbox{$\frac12,\frac12,\frac12,\frac58$} in Theorem~\ref{2=1+1a}]
\[ \frac{288}{\pi(\sqrt2-1)}
=\sum_{k=0}^{\infty}
\hyp{c}
{ \frac{1}{2},\frac{3}{8},\frac{5}{8},\frac{5}{8},\frac{5}{16},\frac{13}{16}}
{\\[-3mm] 1,1,1,\frac{7}{8},\frac{17}{16},\frac{25}{16}}_k
\frac{15360 k^3+12736 k^2+2656 k+195}{16^k}.\]
\end{exam}

%%%%%%%%%%%%%%%%%%%%%%%%%%%%%%%%%%%%%%%%%%%%%%%%%%%%%%%%%%%%%%%%%%%%%%%%%%%
\begin{exam}[\fbox{$\frac12,\frac12,\frac12,\frac78$} in Theorem~\ref{2=1+1a}]
\[ \frac{1056}{\pi(\sqrt2+1)}
=\sum_{k=0}^{\infty}
\hyp{c}
{ \frac{1}{2},\frac{1}{8},\frac{7}{8},\frac{7}{8},\frac{7}{16},\frac{15}{16}}
{\\[-3mm] 1,1,1,\frac{5}{8},\frac{19}{16},\frac{27}{16}}_k
\frac{15360 k^3+14144 k^2+2656 k+105}{16^k}.\]
\end{exam}

%%%%%%%%%%%%%%%%%%%%%%%%%%%%%%%%%%%%%%%%%%%%%%%%%%%%%%%%%%%%%%%%%%%%%%%%%%%
\begin{exam}[\fbox{$\frac12,\frac12,\frac12,\frac16$} in Theorem~\ref{2=1+1b}]
\[ \frac{10 \sqrt{3}}{\pi }
=\sum_{k=0}^{\infty}
\hyp{c}
{ -\frac{1}{2},\frac{1}{6},\frac{1}{6},\frac{1}{12},-\frac{5}{12}}
{\\[-3mm] 1,~ 1,~ 1,~ \frac{1}{3},~\frac{1}{3}}_k
\frac{2160 k^3-372 k^2+68 k+5}{16^k}.\]
\end{exam}

%%%%%%%%%%%%%%%%%%%%%%%%%%%%%%%%%%%%%%%%%%%%%%%%%%%%%%%%%%%%%%%%%%%%%%%%%%%
\begin{exam}[\fbox{$\frac12,\frac12,\frac12,-\frac13$} in Theorem~\ref{2=1+1b}]
\[ \frac{6 \sqrt{3}}{\pi }
=\sum_{k=0}^{\infty}
\hyp{c}
{ -\frac{1}{2},\frac{1}{3},-\frac{1}{3},-\frac{1}{3},-\frac{2}{3}}
{\\[-3mm] 1,~ 1,~ 1,~ \frac{1}{12},~\frac{7}{12}}_k
\frac{135 k^3-48 k^2-7 k+2}{16^k}.\]
\end{exam}

%%%%%%%%%%%%%%%%%%%%%%%%%%%%%%%%%%%%%%%%%%%%%%%%%%%%%%%%%%%%%%%%%%%%%%%%%%%
\begin{exam}[\fbox{$\frac12,\frac12,\frac12,-\frac14$} in Theorem~\ref{2=1+1b}]
\[ \frac{20}{\pi }
=\sum_{k=0}^{\infty}
\hyp{c}
{ -\frac{1}{2},\frac{1}{4},-\frac{1}{4},-\frac{1}{8},-\frac{5}{8}}
{\\[-3mm] 1,~ 1,~ 1,~ \frac{1}{8},~\frac{5}{8}}_k
\frac{960 k^3-232 k^2-38 k+5}{16^k}.\]
\end{exam}

%%%%%%%%%%%%%%%%%%%%%%%%%%%%%%%%%%%%%%%%%%%%%%%%%%%%%%%%%%%%%%%%%%%%%%%%%%%
\begin{exam}[\fbox{$\frac12,\frac12,\frac16,\frac12$} in Theorem~\ref{2=1+1b}]
\[ \frac{10\sqrt{3}}{9  \pi }
=\sum_{k=0}^{\infty}
\hyp{c}
{ \frac{1}{2},\frac{1}{4},-\frac{1}{4},\frac{5}{6},-\frac{5}{6}}
{\\[-3mm] 1,~ 1,~ 1,~ \frac{1}{3},~\frac{2}{3}}_k
\frac{k (120 k^2-26 k+5)}{16^k}.\]
\end{exam}

%%%%%%%%%%%%%%%%%%%%%%%%%%%%%%%%%%%%%%%%%%%%%%%%%%%%%%%%%%%%%%%%%%%%%%%%%%%
\begin{exam}[\fbox{$\frac12,\frac12,\frac76,-\frac12$} in Theorem~\ref{2=1+1b}]
\[ \frac{6 \sqrt{3}}{\pi }
=\sum_{k=0}^{\infty}
\hyp{c}
{ -\frac{1}{2},-\frac{1}{4},-\frac{3}{4},\frac{1}{6},-\frac{1}{6}}
{\\[-3mm] 1,~ 1,~ 1,~ \frac{1}{3},~\frac{2}{3}}_k
\frac{720 k^3-300 k^2-4 k+3}{16^k}.\]
\end{exam}

%%%%%%%%%%%%%%%%%%%%%%%%%%%%%%%%%%%%%%%%%%%%%%%%%%%%%%%%%%%%%%%%%%%%%%%%%%%
\begin{exam}[\fbox{$\frac12,\frac12,-\frac12,\frac76$} in Theorem~\ref{2=1+1b}]
\[ \frac{27 \sqrt{3}}{\pi }
=\sum_{k=0}^{\infty}
\hyp{c}
{ -\frac{3}{2},\frac{1}{6},\frac{7}{6},\frac{1}{12},\frac{7}{12}}
{\\[-3mm] 1,~ 1,~ 1,~ \frac{1}{3},-\frac{2}{3}}_k
\frac{21+292 k-420 k^2-2160 k^3}{16^k}.\]
\end{exam}

%%%%%%%%%%%%%%%%%%%%%%%%%%%%%%%%%%%%%%%%%%%%%%%%%%%%%%%%%%%%%%%%%%%%%%%%%%%
\begin{exam}[\fbox{$\frac12,\frac12,-\frac12,\frac14$} in Theorem~\ref{2=1+1b}]
\[ \frac{96}{\pi }
=\sum_{k=0}^{\infty}
\hyp{c}
{ -\frac{3}{2},\frac{1}{4},\frac{3}{4},\frac{1}{8},-\frac{3}{8}}
{\\[-3mm] 1,~ 1,~ 1,~\frac{3}{8},-\frac{1}{8}}_k
\frac{9+102 k-424 k^2-960 k^3}{16^k}.\]
\end{exam}

%%%%%%%%%%%%%%%%%%%%%%%%%%%%%%%%%%%%%%%%%%%%%%%%%%%%%%%%%%%%%%%%%%%%%%%%%%%
\begin{exam}[\fbox{$\frac12,\frac12,-\frac12,-\frac14$} in Theorem~\ref{2=1+1b}]
\[ \frac{160}{\pi }
=\sum_{k=0}^{\infty}
\hyp{c}
{ -\frac{3}{2},\frac{5}{4},-\frac{5}{4},-\frac{1}{8},-\frac{5}{8}}
{\\[-3mm] 1,~ 1,~ 1,~ \frac{1}{8},-\frac{3}{8}}_k
\frac{960 k^3-232 k^2-710 k+75 }{16^k}.\]
\end{exam}

%%%%%%%%%%%%%%%%%%%%%%%%%%%%%%%%%%%%%%%%%%%%%%%%%%%%%%%%%%%%%%%%%%%%%%%%%%%
\begin{exam}[\fbox{$\frac12,\frac12,-\frac16,\frac12$} in Theorem~\ref{2=1+1b}]
\[ \frac{28}{3 \sqrt{3} \pi }
=\sum_{k=0}^{\infty}
\hyp{c}
{ \frac{1}{2},\frac{1}{4},-\frac{1}{4},\frac{7}{6},-\frac{7}{6},\frac{13}{12}}
{\\[-3mm] 1,1,1,~ \frac{2}{3},\frac{4}{3},~\frac{1}{12}}_k
\frac{k (60 k^2-8 k-7)}{16^k}.\]
\end{exam}

%%%%%%%%%%%%%%%%%%%%%%%%%%%%%%%%%%%%%%%%%%%%%%%%%%%%%%%%%%%%%%%%%%%%%%%%%%%
\begin{exam}[\fbox{$\frac12,\frac12,-\frac12,-\frac13$} in Theorem~\ref{2=1+1b}]
\[ \frac{162 \sqrt{3}}{5 \pi }
=\sum_{k=0}^{\infty}
\hyp{c}
{ -\frac{3}{2},-\frac{1}{3},-\frac{2}{3},-\frac{4}{3},\frac{4}{3}}
{\\[-3mm] 1,1,1,~ \frac{1}{12},-\frac{5}{12}}_k
\frac{135 k^3-48 k^2-106 k+24}{16^k}.\]
\end{exam}

\subsection{Series for $\pi$}
%%%%%%%%%%%%%%%%%%%%%%%%%%%%%%%%%%%%%%%%%%%%%%%%%%%%%%%%%%%%%%%%%%%%%%%%%%%
%%%%%%%%%%%%%%%%%%%%%%%%%%%%%%%%%%%%%%%%%%%%%%%%%%%%%%%%%%%%%%%%%%%%%%%%%%%
%%%%%%%%%%%%%%%%%%%%%%%%%%%%%%%%%%%%%%%%%%%%%%%%%%%%%%%%%%%%%%%%%%%%%%%%%%%

%%%%%%%%%%%%%%%%%%%%%%%%%%%%%%%%%%%%%%%%%%%%%%%%%%%%%%%%%%%%%%%%%%%%%%%%%%%
\begin{exam}[\fbox{$\frac52,2,2,\frac34$} in Theorem~\ref{2=1+1b}]
\[ \frac{5\pi}{16}
=\sum_{k=0}^{\infty}\hyp{c}
{-\frac{1}{2},-\frac{1}{4},\frac{1}{8},-\frac{3}{8}}
{\\[-3mm] \frac{1}{2},~\frac{3}{4},~\frac{9}{8},~\frac{13}{8}}_k
\frac{1-7k+40k^2}{16^k}.\]
\end{exam}

%%%%%%%%%%%%%%%%%%%%%%%%%%%%%%%%%%%%%%%%%%%%%%%%%%%%%%%%%%%%%%%%%%%%%%%%%%%
\begin{exam}[\fbox{$\frac52,2,2,\frac14$} in Theorem~\ref{2=1+1b}]
\[ \frac{25 \pi }{16}
=\sum_{k=0}^{\infty}\hyp{c}
{ -\frac{1}{2},\frac{1}{4},-\frac{1}{8},-\frac{5}{8}}
{\\[-3mm] \frac{1}{2},~\frac{9}{4},~\frac{7}{8},~\frac{11}{8}}_k
\frac{120 k^2+77 k+5}{16^k}.\]
\end{exam}

%%%%%%%%%%%%%%%%%%%%%%%%%%%%%%%%%%%%%%%%%%%%%%%%%%%%%%%%%%%%%%%%%%%%%%%%%%%
\begin{exam}[\fbox{$\frac32,2,1,\frac14$} in Theorem~\ref{2=1+1b}]
\[ \frac{3\pi}{8}
=\sum_{k=0}^{\infty}\hyp{c}
{-\frac{1}{2},\frac{1}{4},-\frac{1}{8},\frac{3}{8}}
{\\[-3mm] \frac{1}{2},\:\frac{5}{4},\:\frac{7}{8},\frac{11}{8}}_k
\frac{1+11 k+106 k^2+240 k^3}{16^k}.\]
\end{exam}

\begin{exam}[\fbox{$\frac32,1,1,\frac56$} in Theorem~\ref{2=1+1a}]
\[ \frac{36 \pi }{5 \sqrt{3}}
=\sum_{k=0}^{\infty}
\hyp{c}
{ \frac{1}{2},\frac{1}{3},\frac{2}{3},\frac{1}{6},\frac{5}{6},\frac{14}{9}}
{\\[-3mm] \frac{3}{2},\frac{5}{3},\frac{5}{4},\frac{7}{4},\frac{7}{6},\frac{5}{9}}_k
\frac{60 k^2+64 k+13}{16^k}.\]
\end{exam}

%%%%%%%%%%%%%%%%%%%%%%%%%%%%%%%%%%%%%%%%%%%%%%%%%%%%%%%%%%%%%%%%%%%%%%%%%%%%
%\begin{exam}[\fbox{\bf REPLACE}\fbox{$\frac12,2,1,-\frac16$} in Theorem~\ref{2=1+1a}]
%\[-\frac{140 \pi }{297 \sqrt{3}}
%=\sum_{k=0}^{\infty}
%\hyp{c}
%{ 2,\frac{5}{2},-\frac{3}{2},\frac{1}{3},\frac{2}{3},-\frac{1}{6},\frac{7}{6},\frac{41}{30}}
%{\\[-3mm]  1,-\frac{1}{2},\frac{3}{2},\frac{5}{3},\frac{3}{4},\frac{5}{4},\frac{13}{6},\frac{11}{30}}_k
%\frac{1+38k^2+36k^3}{16^k}.\]
%\end{exam}
%
%%%%%%%%%%%%%%%%%%%%%%%%%%%%%%%%%%%%%%%%%%%%%%%%%%%%%%%%%%%%%%%%%%%%%%%%%%%%
%\begin{exam}[\fbox{\bf REPLACE}\fbox{$\frac52,1,2,-\frac23$} in Theorem~\ref{2=1+1a}]
%\[\frac{513\sqrt3\pi }{6400}
%=\sum_{k=0}^{\infty}
%\hyp{c}
%{ \frac{1}{2},\frac{8}{3},-\frac{2}{3},-\frac{7}{6},\frac{13}{6}}
%{\\[-3mm] \frac{5}{2},~\frac{5}{3},~\frac{5}{4},~\frac{7}{4},~\frac{25}{6}}_k
%\frac{50+179 k-834 k^2-1080 k^3}{16^k(1-12k)(7-12k)(13-12k)}.\]
%
%\[\frac{567\pi }{6400\sqrt{3}}
%=\sum_{k=0}^{\infty}
%\hyp{c}
%{ \frac{1}{2},\frac{8}{3},-\frac{2}{3},-\frac{7}{6},-\frac{7}{12},-\frac{13}{12}}
%{\\[-3mm] \frac{5}{2},~\frac{5}{3},~\frac{5}{4},~\frac{7}{4},~\frac{5}{12},~\frac{11}{12}}_k
%\frac{50+179 k-834 k^2-1080 k^3}{16^k(6k+13)(6k+19)}.\]
%\end{exam}

%%%%%%%%%%%%%%%%%%%%%%%%%%%%%%%%%%%%%%%%%%%%%%%%%%%%%%%%%%%%%%%%%%%%%%%%%%%
\begin{exam}[Chu and Zhang~\cito{chu14mc}:
\fbox{$\frac32,1,\frac56,1$} in Theorem~\ref{2=1+1a}]
\[ \frac{20 \pi }{9 \sqrt{3}}
=\sum_{k=0}^{\infty}
\hyp{c}
{ 1,\frac{1}{3},\frac{2}{3},~\frac{1}{4},\frac{3}{4},\frac{8}{5}}
{\\[-3mm] \frac{3}{2},\frac{3}{2},\frac{3}{2},\frac{3}{5},\frac{7}{6},\frac{11}{6}}_k
\frac{12 k^2+15 k+4}{16^k}.\]
\end{exam}

%%%%%%%%%%%%%%%%%%%%%%%%%%%%%%%%%%%%%%%%%%%%%%%%%%%%%%%%%%%%%%%%%%%%%%%%%%%
\begin{exam}[\fbox{$\frac12,1,\frac16,1$} in Theorem~\ref{2=1+1a}]
\[\frac{20 \pi }{27 \sqrt{3}}
=\sum_{k=0}^{\infty}
\hyp{c}
{ 1,\frac{1}{3},\frac{2}{3},~\frac{3}{4},\frac{5}{4}}
{\\[-3mm] \frac{1}{2},\frac{1}{2},\frac{3}{2},\frac{7}{6},\frac{11}{6}}_k
\frac{4-11 k- 69 k^2 -60 k^3}{16^k}.\]
\end{exam}

%%%%%%%%%%%%%%%%%%%%%%%%%%%%%%%%%%%%%%%%%%%%%%%%%%%%%%%%%%%%%%%%%%%%%%%%%%%
\begin{exam}[\fbox{$\frac32,1,\frac56,2$} in Theorem~\ref{2=1+1a}]
\[ \frac{140 \pi }{27 \sqrt{3}}
=\sum_{k=0}^{\infty}
\hyp{c}
{ 2,\frac{4}{3},-\frac{1}{3},\frac{3}{4},\frac{5}{4}}
{\\[-3mm] \frac{1}{2},\frac{3}{2},\frac{3}{2},\frac{11}{6},\frac{13}{6}}_k
\frac{60 k^3+133 k^2+85 k+13}{16^k}.\]
\end{exam}

%%%%%%%%%%%%%%%%%%%%%%%%%%%%%%%%%%%%%%%%%%%%%%%%%%%%%%%%%%%%%%%%%%%%%%%%%%%
\begin{exam}[\fbox{$\frac12,1,-\frac16,2$} in Theorem~\ref{2=1+1a}]
\[\frac{700 \pi }{243 \sqrt{3}}
=\sum_{k=0}^{\infty}
\hyp{c}
{  2,-\frac{1}{3},\frac{4}{3},~\frac{5}{4},\frac{7}{4}}
{\\[-3mm]  -\frac{1}{2},\frac{1}{2},\frac{3}{2},\frac{11}{6},\frac{13}{6}}_k
\frac{25-3 k-91 k^2-60 k^3}{16^k}.\]
\end{exam}

%%%%%%%%%%%%%%%%%%%%%%%%%%%%%%%%%%%%%%%%%%%%%%%%%%%%%%%%%%%%%%%%%%%%%%%%%%%
\begin{exam}[\fbox{$\frac32,1,\frac56,1$} in Theorem~\ref{2=1+1b}]
\[\frac{4 \pi }{9 \sqrt{3}}
=\frac23+\sum_{k=1}^{\infty}
\hyp{c}
{ 1,\frac{2}{3},-\frac{2}{3},\frac{1}{4},-\frac{1}{4} }
{\\[-3mm] \frac{1}{2},\frac{1}{2},\frac{3}{2},~\frac{5}{6},\frac{7}{6}}_k
\frac{20 k^2+7 k+2}{16^k}.\]
\end{exam}

%%%%%%%%%%%%%%%%%%%%%%%%%%%%%%%%%%%%%%%%%%%%%%%%%%%%%%%%%%%%%%%%%%%%%%%%%%%
\begin{exam}[\fbox{$\frac12,1,\frac16,1$} in Theorem~\ref{2=1+1b}]
\[\frac{4 \pi }{81 \sqrt{3}}
=\frac23+\sum_{k=1}^{\infty}
\hyp{c}
{ 1,~ \frac{2}{3},-\frac{2}{3},\frac{1}{4},\frac{3}{4}}
{\\[-3mm]\frac{1}{2},\frac{1}{2},-\frac{1}{2},\frac{5}{6},\frac{7}{6}}_k
\frac{2+7k-20 k^2}{16^k}.\]
\end{exam}

%%%%%%%%%%%%%%%%%%%%%%%%%%%%%%%%%%%%%%%%%%%%%%%%%%%%%%%%%%%%%%%%%%%%%%%%%%%
\begin{exam}[\fbox{$\frac52,1,\frac56,2$} in Theorem~\ref{2=1+1b}]
\[ \frac{10 \pi }{7 \sqrt{3}}
=2+\sum_{k=1}^{\infty}
\hyp{c}
{2,\frac{5}{3},-\frac{5}{3},\frac{1}{4},-\frac{1}{4}}
{\\[-3mm] \frac{1}{2},\frac{1}{2},\frac{5}{2},~\frac{7}{6},\frac{11}{6} }_k
\frac{12 k^2+17 k+8}{16^k}.\]
\end{exam}

%%%%%%%%%%%%%%%%%%%%%%%%%%%%%%%%%%%%%%%%%%%%%%%%%%%%%%%%%%%%%%%%%%%%%%%%%%%
\begin{exam}[\fbox{$\frac32,1,2,\frac56$} in Theorem~\ref{2=1+1b}]
\[\frac{16 \pi }{3 \sqrt{3}}
=20+\sum_{k=1}^{\infty}
\hyp{c}
{ -\frac{1}{2},\frac{4}{3},-\frac{4}{3},\frac{1}{6},\frac{5}{6}}
{\\[-3mm] \frac{1}{2},\frac{5}{3},~\frac{1}{4},\frac{3}{4},\frac{7}{6}}_k
\frac{56-33 k-270 k^2}{16^k}.\]
\end{exam}

%%%%%%%%%%%%%%%%%%%%%%%%%%%%%%%%%%%%%%%%%%%%%%%%%%%%%%%%%%%%%%%%%%%%%%%%%%%
\begin{exam}[\fbox{$\frac12,1,\frac16,2$} in Theorem~\ref{2=1+1b}]
\[\frac{350 \pi }{243 \sqrt{3}}
=30+\sum_{k=1}^{\infty}
\hyp{c}
{2,~\frac{5}{3},-\frac{5}{3},\frac{3}{4},\frac{5}{4}}
{\\[-3mm]\frac{1}{2},\frac{1}{2},-\frac{3}{2},\frac{7}{6},\frac{11}{6}}_k
\frac{40-k-60 k^2}{16^k}.\]
\end{exam}

%%%%%%%%%%%%%%%%%%%%%%%%%%%%%%%%%%%%%%%%%%%%%%%%%%%%%%%%%%%%%%%%%%%%%%%%%%%
\begin{exam}[\fbox{$\frac32,1,\frac16,2$} in Theorem~\ref{2=1+1b}]
\[\frac{32 \pi }{27 \sqrt{3}}
=-8+\sum_{k=1}^{\infty}
\hyp{c}
{ 2,-\frac{2}{3},-\frac{4}{3},\frac{1}{4},\frac{3}{4}}
{\\[-3mm] \frac{1}{2},-\frac{1}{2},\frac{3}{2},~\frac{5}{6},\frac{7}{6}}_k
\frac{k(15k+2)(1-12k)}{16^k}.\]
\end{exam}

%%%%%%%%%%%%%%%%%%%%%%%%%%%%%%%%%%%%%%%%%%%%%%%%%%%%%%%%%%%%%%%%%%%%%%%%%%%
\begin{exam}[\fbox{$\frac52,1,2,\frac56$} in Theorem~\ref{2=1+1a}]
\[ \frac{270\pi }{7\sqrt{3}}
=\sum_{k=0}^{\infty}
\hyp{c}
{ \frac{1}{2},\frac{1}{3},-\frac{1}{3},\frac{1}{6},\frac{5}{6}}
{\\[-3mm] \frac{5}{2},\frac{8}{3},~\frac{5}{4},\frac{7}{4},\frac{7}{6}}_k
\frac{70+409 k+627 k^2+270 k^3}{16^k}.\]
\end{exam}

%%%%%%%%%%%%%%%%%%%%%%%%%%%%%%%%%%%%%%%%%%%%%%%%%%%%%%%%%%%%%%%%%%%%%%%%%%%
\begin{exam}[\fbox{$\frac52,1,2,\frac16$} in Theorem~\ref{2=1+1a}]
\[ \frac{756\pi }{275\sqrt{3}}
=\sum_{k=0}^{\infty}
\hyp{c}
{ \frac{1}{2},-\frac{1}{3},-\frac{2}{3},\frac{1}{6},-\frac{1}{6},\frac{11}{6}}
{\\[-3mm] \frac{5}{2},~\frac{10}{3},~\frac{5}{4},~\frac{7}{4},~\frac{5}{6},~\frac{5}{6}}_k
\frac{5+92 k+258 k^2+135 k^3}{16^k}.\]
\end{exam}

\subsection{\textbf{BBP}--series} \
%%%%%%%%%%%%%%%%%%%%%%%%%%%%%%%%%%%%%%%%%%%%%%%%%%%%%%%%%%%%%%%%%%%%%%%%%%%
%%%%%%%%%%%%%%%%%%%%%%%%%%%%%%%%%%%%%%%%%%%%%%%%%%%%%%%%%%%%%%%%%%%%%%%%%%%
%%%%%%%%%%%%%%%%%%%%%%%%%%%%%%%%%%%%%%%%%%%%%%%%%%%%%%%%%%%%%%%%%%%%%%%%%%%
In 1995, Simon Plouffe discovered the following amazing BBP--formula 
(named after Bailey--Borwein--Plouffe~\citu{bbp97mc}{Theorem~1})
\[\pi=\sum_{k=0}^{\infty}\Big(\frac1{16}\Big)^k
\bigg\{\frac{4}{8 k+1}-\frac{2}{8 k+4}
-\frac{1}{8 k+5}-\frac{1}{8 k+6}\bigg\}\]
that provides a digit-extraction algorithm for $\pi$ in base 10. 
By decomposing the factorial fraction in the summand into partial 
fractions, we can show that the next five series are all equivalent 
to the above BBP--formula.

%%%%%%%%%%%%%%%%%%%%%%%%%%%%%%%%%%%%%%%%%%%%%%%%%%%%%%%%%%%%%%%%%%%%%%%%%%%
\begin{exam}[\fbox{$\frac32,1,1,\frac34$} in Theorem~\ref{2=1+1a}]
\[ 15 \pi
=\sum_{k=0}^{\infty}\hyp{c}
{\frac{1}{2},\frac{3}{4},\frac{1}{8},\:\frac{5}{8}}
{\\[-3mm] \frac{3}{2},\frac{7}{4},\frac{9}{8},\frac{13}{8}}_k
\frac{120 k^2+151 k+47}{16^k}.\]
\end{exam}

%%%%%%%%%%%%%%%%%%%%%%%%%%%%%%%%%%%%%%%%%%%%%%%%%%%%%%%%%%%%%%%%%%%%%%%%%%%
\begin{exam}[\fbox{$\frac52,1,2,\frac34$} in Theorem~\ref{2=1+1a}]
\[ \frac{63 \pi }{2}
=\sum_{k=0}^{\infty}\hyp{c}
{ \frac{1}{2},\frac{3}{4},\frac{1}{8},-\frac{3}{8}}
{\\[-3mm] \frac{5}{2},\frac{11}{4},\frac{9}{8},\frac{13}{8}}_k
\frac{120 k^2+235 k+99 }{16^k}.\]
\end{exam}

%%%%%%%%%%%%%%%%%%%%%%%%%%%%%%%%%%%%%%%%%%%%%%%%%%%%%%%%%%%%%%%%%%%%%%%%%%%
\begin{exam}[\fbox{$\frac32,1,2,-\frac14$} in Theorem~\ref{2=1+1b}]
\[ \frac{21\pi}{8}
=7+\sum_{k=1}^{\infty}\hyp{c}
{-\frac{1}{2},-\frac{1}{4},-\frac{3}{8},-\frac{7}{8}}
{\\[-3mm] \frac{1}{2},~\frac{7}{4},~\frac{5}{8},~\frac{9}{8}}_k
\frac{480k^2-172k-9}{16^k}.\]
\end{exam}

%%%%%%%%%%%%%%%%%%%%%%%%%%%%%%%%%%%%%%%%%%%%%%%%%%%%%%%%%%%%%%%%%%%%%%%%%%%
\begin{exam}[\fbox{$\frac52,1,2,\frac34$} in Theorem~\ref{2=1+1b}]
\[\frac{21 \pi }{10}
=7+\sum_{k=1}^{\infty}\hyp{c}
{-\frac{1}{2},-\frac{1}{4},-\frac{3}{8},-\frac{7}{8}}
{\\[-3mm] \frac{5}{2},~\frac{3}{4},~\frac{5}{8},~\frac{9}{8} }_k
\frac{23+10 k-240 k^2}{16^k}.\]
\end{exam}

%%%%%%%%%%%%%%%%%%%%%%%%%%%%%%%%%%%%%%%%%%%%%%%%%%%%%%%%%%%%%%%%%%%%%%%%%%%
\begin{exam}[\fbox{$\frac32,1,2,-\frac54$} in Theorem~\ref{2=1+1b}]
\[ \frac{77\pi}{8}
=-\frac{55}{3}+\sum_{k=1}^{\infty}\hyp{c}
{-\frac{1}{2},-\frac{5}{4},-\frac{7}{8},-\frac{11}{8}}
{\\[-3mm] \frac{1}{2},~\frac{11}{4},~\frac{1}{8},~\frac{5}{8}}_k
\frac{160k^2-36k-13}{16^k}.\]
\end{exam}

%%%%%%%%%%%%%%%%%%%%%%%%%%%%%%%%%%%%%%%%%%%%%%%%%%%%%%%%%%%%%%%%%%%%%%%%%%%
%%%%%%%%%%%%%%%%%%%%%%%%%%%%%%%%%%%%%%%%%%%%%%%%%%%%%%%%%%%%%%%%%%%%%%%%%%%
%%%%%%%%%%%%%%%%%%%%%%%%%%%%%%%%%%%%%%%%%%%%%%%%%%%%%%%%%%%%%%%%%%%%%%%%%%%
There is another BBP--formula disguised in the article 
by Adamchik--Wagon~\cito{adam+w}
\[2\pi=\sum_{k=0}^{\infty}\Big(\frac1{16}\Big)^k
\bigg\{\frac{8}{8 k+2}+\frac{4}{8 k+3}+\frac{4}{8 k+4}-\frac{1}{8 k+7}\bigg\}.\]
Then the same approach of partial fractions can show that it 
has the following different infinite series representations.

%%%%%%%%%%%%%%%%%%%%%%%%%%%%%%%%%%%%%%%%%%%%%%%%%%%%%%%%%%%%%%%%%%%%%%%%%%%
\begin{exam}[\fbox{$\frac32,1,1,\frac14$} in Theorem~\ref{2=1+1b}]
\[\frac{5 \pi }{9}
=\frac53+\sum_{k=1}^{\infty}\hyp{c}
{ -\frac{1}{2},\frac{1}{4},-\frac{1}{8},-\frac{5}{8}}
{\\[-3mm] \frac{3}{2},~\frac{5}{4},~\frac{3}{8},~\frac{7}{8}}_k
\frac{7-6k-80 k^2}{16^k}.\]
\end{exam}

%%%%%%%%%%%%%%%%%%%%%%%%%%%%%%%%%%%%%%%%%%%%%%%%%%%%%%%%%%%%%%%%%%%%%%%%%%%
\begin{exam}[\fbox{$\frac52,1,2,\frac14$} in Theorem~\ref{2=1+1b}]
\[\frac{15 \pi }{14}
=3+\sum_{k=1}^{\infty}\hyp{c}
{ -\frac{1}{2},\frac{1}{4},-\frac{5}{8},-\frac{9}{8}}
{\\[-3mm] \frac{5}{2},~\frac{9}{4},~\frac{3}{8},~\frac{7}{8} }_k
\frac{19-62 k-80 k^2}{16^k}.\]
\end{exam}

\begin{exam}[\fbox{$\frac32,1,2,\frac14$} in Theorem~\ref{2=1+1b}]
\[\frac{15 \pi }{8}
=5+\sum_{k=1}^{\infty}\hyp{c}
{-\frac{1}{2},-\frac{3}{4},-\frac{1}{8},-\frac{5}{8}}
{\\[-3mm] \frac{1}{2},~\frac{1}{4},~\frac{7}{8},~\frac{11}{8}}_k
\frac{160 k^2-108 k+21}{16^k}.\]
\end{exam}

%%%%%%%%%%%%%%%%%%%%%%%%%%%%%%%%%%%%%%%%%%%%%%%%%%%%%%%%%%%%%%%%%%%%%%%%%%%
\begin{exam}[\fbox{$\frac32,1,2,-\frac34$} in Theorem~\ref{2=1+1b}]
\[ \frac{45\pi}{8}
=16+\sum_{k=0}^{\infty}\hyp{c}
{-\frac{1}{2},-\frac{3}{4},-\frac{5}{8},-\frac{9}{8}}
{\\[-3mm] \frac{1}{2},~\frac{9}{4},~\frac{3}{8},~\frac{7}{8}}_k
\frac{11+260k-480k^2}{16^k}.\]
\end{exam}

%%%%%%%%%%%%%%%%%%%%%%%%%%%%%%%%%%%%%%%%%%%%%%%%%%%%%%%%%%%%%%%%%%%%%%%%%%%
\begin{exam}[\fbox{$\frac32,1,1,\frac54$} in Theorem~\ref{2=1+1a}]
\[{21 \pi}=\sum_{k=0}^{\infty}
\hyp{c}{\frac{1}{2},\:\frac{1}{4},\:\frac{3}{8},\:\frac{7}{8}}
{\\[-3mm] \frac{3}{2},\frac{5}{4},\frac{11}{8},\frac{15}{8}}_k
\frac{65+413 k+812 k^2+480 k^3}{16^k}.\]
\end{exam}

%%%%%%%%%%%%%%%%%%%%%%%%%%%%%%%%%%%%%%%%%%%%%%%%%%%%%%%%%%%%%%%%%
%%%%%%%%%%%%%%%%%%%%%%%%%%%%%%%%%%%%%%%%%%%%%%%%%%%%%%%%%%%%%%%%%
%%%%%%%%%%%%%%%%%%%%%%%%%%%%%%%%%%%%%%%%%%%%%%%%%%%%%%%%%%%%%%%%%

%%%%%%%%%%%%%%%%%%%%%%%%%%%%%%%%%%%%%%%%%%%%%%%%%%%%%%%%%%%%%%
\end{document}